\documentclass[leqno]{article}
\usepackage{amssymb}
\usepackage{amsmath, enumerate, vmargin}
\usepackage[arrow, matrix]{xy}

\makeatletter
\renewcommand\section{\@startsection {section}{1}{\z@}%
 {-3.5ex \@plus -1ex \@minus -.2ex}%
 {2.3ex \@plus.2ex}%
 {\center \normalfont\large\bfseries}}
\makeatother

\setmarginsrb{3cm}{3cm}{3cm}{3cm}{1mm}{6mm}{0mm}{10mm}

\newtheorem{thm}{Theorem}[section]
\newtheorem{prop}[thm]{Proposition}
\newtheorem{cor}[thm]{Corollary}
\newtheorem{lem}[thm]{Lemma}
\newtheorem{defi}[thm]{Definition}
\newtheorem{remark}[thm]{Remark}
\newtheorem{example}[thm]{Example}
\newtheorem{pb}[thm]{Problem}

\newenvironment{rk}{\begin{remark}\rm}{\end{remark}}
\newenvironment{definition}{\begin{defi}\rm}{\end{defi}}

\newcommand{\rz}{{\mathbb R}}
\newcommand{\nz}{{\mathbb N}}
\newcommand{\zz}{{\mathbb Z}}
\newcommand{\cz}{{\mathbb C}}
\newcommand{\un}{1\mkern -4mu{\textrm l}}

\newcommand{\F}{\Phi}

\newcommand{\M}{{\mathcal M}}
\newcommand{\N}{{\mathcal N}}

\newcommand{\U}{{\mathcal U}}

\renewcommand{\a}{\alpha}
\renewcommand{\b}{\beta}
\newcommand{\g}{\gamma}
\newcommand{\Ga}{\Gamma}
\renewcommand{\d}{\delta}

\newcommand{\e}{\varepsilon}

\renewcommand{\l}{\lambda}
\newcommand{\La}{\Lambda}
\newcommand{\f}{\varphi}
\renewcommand{\O}{{\Omega}}
\renewcommand{\o}{{\omega}}
\newcommand{\s}{\sigma}

\newcommand{\Tr}{{\rm Tr}}

\newcommand{\op}{\oplus}
\newcommand{\ot}{\otimes}
\newcommand{\el}{\ell}

\newcommand{\wt}{\widetilde}
\newcommand{\wh}{\widehat}

\newcommand{\id}{{\rm id}}

\newcommand{\xspace}{\hbox{\kern-2.5pt}}

\newcommand{\yspace}{\hbox{\kern-0.9pt}}

\newcommand{\tnorm}[1]{\left\vert\xspace\left\vert\xspace\left\vert\mskip2mu
#1\mskip2mu \right\vert\xspace\right\vert\xspace\right\vert}

\newcommand{\pf}{\noindent{\it Proof.~~}}
\newcommand{\cqd}{\hfill$\Box$}
\newcommand{\be}{\begin{eqnarray*}}
\newcommand{\ee}{\end{eqnarray*}}
\newcommand{\beq}{\begin{equation}}
\newcommand{\eeq}{\end{equation}}
\newcommand{\8}{\infty}
\newcommand{\n}{\noindent}

\numberwithin{equation}{section}

\begin{document}

\title{Representation of certain homogeneous Hilbertian
operator spaces and applications}
\author{Marius Junge\footnotemark\;\; and Quanhua Xu\footnotemark}
\date{}

\maketitle

\vskip 1.5cm
\begin{abstract}
Following Grothendieck's characterization of Hilbert spaces we
consider operator spaces $F$ such that both $F$ and $F^*$
completely embed into the dual of a C*-algebra. Due to
Haagerup/Musat's improved version  of Pisier/Shlyakhtenko's
Grothendieck inequality  for operator spaces, these spaces are
quotients of subspaces of the direct sum $C\oplus R$ of the column
and row spaces (the corresponding class being denoted by $QS(C\op
R)$). We first prove a representation theorem for homogeneous
$F\in QS(C\op R)$ starting from the fundamental sequences
  $$\F_{c}(n)=
 \big\|\sum_{k=1}^n e_{k 1}\ot e_k\big\|_{C\ot_{\min}F}^2
 \quad\mbox{and}
 \quad \F_{r}(n)=
 \big\|\sum_{k=1}^n e_{1k}\ot e_k\big\|_{R\ot_{\min}F}^2\,$$
given by an orthonormal basis $(e_k)$ of $F$. Under a mild
regularity assumption on these sequences we show that they
completely determine the operator space structure of $F$ and find
a canonical representation of this important class of homogeneous
Hilbertian operator spaces in terms of weighted row and column
spaces. This canonical representation allows us to get an explicit
formula for the exactness constant of an $n$-dimensional subspace
$F_n$ of $F$:
 $$ex(F_n)\sim \left[\frac{n}{\F_{c}(n)}\,
  \F_{r}\big(\frac{\F_{c}(n)}{\F_{r}(n)}\big) +
  \frac{n}{\F_{r}(n)}\,
  \F_{c}\big(\frac{\F_{r}(n)}{\F_{c}(n)}\big)
  \right]^{1/2}.$$
In the same way, the projection (=injectivity) constant of $F_n$
is explicitly expressed in terms of $\Phi_c$ and $\Phi_r$ too.
Orlicz space techniques play a crucial role in our arguments. They
also permit us to determine the completely 1-summing maps in
Effros and Ruan's sense between two homogeneous spaces $E$ and $F$
in $QS(C\op R)$. The resulting space $\Pi_1^o(E,\, F)$
isomorphically coincides with a Schatten-Orlicz class $S_\f$.
Moreover, the underlying Orlicz function $\f$ is uniquely
determined by the fundamental sequences of $E$ and $F$. In
particular,  applying these results to the column subspace $C_p$
of the Schatten $p$-class, we find the projection and exactness
constants of $C_p^n$,  and determine the completely 1-summing maps
from $C_p$ to $C_q$ for any $1\le p,\, q\le\8$.

 \end{abstract}

 \setcounter{section}{-1}


 \makeatletter
 \renewcommand{\@makefntext}[1]{#1}
 \makeatother \footnotetext{\noindent
 $^*$ Partially supported by NSF DMS 05-56120.\\
  $^\dagger$  Partially supported by ANR 06-BLAN-0015.
 }

\newpage


\section{Introduction}


Hilbert spaces are by far the most important examples among
general Banach spaces. Indeed, much research in Banach space
theory relies on weakening abstract properties of Hilbert spaces
and studies  the class of Banach spaces satisfying these weaker
properties. This applies in particular to the notions of type and
cotype (see \cite{pis-fact}), as well as to spaces with
unconditional martingale differences (so called UMD spaces) or
equivalently to spaces for which the Hilbert transform
continuously extends to the vector-valued setting (see
\cite{burk-UMD}, \cite{burk-HT} and \cite{bourg-UMD})

In the theory of quantized Banach spaces, i.e., the theory of
operator spaces, the class of Hilbert spaces is also quantized. In
other words, for a given Hilbert space $F$ there exist numerous
isometric embeddings of $F$ into the space $B(H)$ of bounded
operators on a Hilbert space $H$, which give rise to very
different matrix norms on $F$. First examples that one should have
in mind are, of course, the column and row spaces
 $$ C=\overline{\rm span}\{e_{k1}\,:\, k\in \nz\}\subset B(\el_2)
  \quad \mbox{and}\quad
 R=\overline{\rm span}\{e_{1k}\,:\, k\in \nz\}\subset B(\el_2)\,.$$
We may think of $C$ as a Hilbert space with a column spin and $R$
with a row spin. Both $C$ and $R$ are isometric to $\el_2$ as
Banach spaces. They are, however, extremely different as operator
spaces. The spaces $C$ and $R$ are omnipresent in
operator space theory. For instance, Pisier's operator Hilbert
space $OH$ can be constructed from them via complex interpolation.

The spaces $C$, $R$ and $OH$ are examples of homogeneous
Hilbertian operator spaces. Recall that an operator space
$F\subset B(H)$ is called \emph{Hilbertian}, respectively,
\emph{homogeneous} if it is isomorphic to a Hilbert space at the
Banach space level, respectively, if every bounded linear map $T$
on $F$ is automatically completely bounded. Let us  point out that
these spaces are dramatically different from the operator space
structures on $\el_2$ constructed by Oikhberg and Ricard
\cite{oik-ric} which allow only very few completely bounded maps.

The importance of $C$ and $R$ in operator space theory is also
illustrated by the following noncommutative analogue of
Grothendieck's abstract characterization of Hilbert spaces in the
commutative theory. For a (separable) operator space  $F$  both
$F$ and $F^*$ completely embed in a noncommutative $L_1$-space if
and only if $F$ is completely isomorphic to a quotient of a
subspace of the direct sum $C\op R$. The ``only if'' part is due
to Pisier/Shlyakhtenko \cite{pisshlyak} with an exactness
assumption and Haagerup/Musat \cite{haag-mu-er} in the full
generality; the ``if'' part was proved independently by the
present authors (see \cite{ju-araki} and \cite{xu-embed}). Let
$QS(C\op R)$ denote the class of quotients of subspaces of $C\op
R$ and $HQS(C\op R)$ the subclass of homogeneous spaces. This last
class is the main concern of the present paper. We will study
various properties of it.

In the literature we find a very particular construction of spaces
in $HQS(C\op R)$ closely related to quasi-free states. Let $u_c$
and $u_r$ be two weights on a measure space $(\Omega, \nu)$. Then
we may consider the weighted column and row spaces
$L_2^c(u_c)=B(\cz,\,L_2(u_c))$, $L_2^r(u_r)=B(L_2(u_r),\,\cz)$ and
their sum
 $$L_2^c(u_c)+L_2^r(u_r)  = \left\{a+b\,:\, a\in L_2^c(u_c),\;
 b\in L_2^r(u_r)\right\}.$$
It turns out that this sum (which is no longer direct) can be
realized in the predual of a quasi-free von Neumann algebra
(either in the free or classical sense). Clearly,
$L_2^c(u_c)+L_2^r(u_r)$ is a quotient of $L_2^c(u_c)\op
L_2^r(u_r)$. However, we may also consider the subspace of
constant functions of $ L_2^c(u_c;\ell_2)+L_2^r(u_r;\ell_2)$:
 $$K_{u_c,\,u_r} = \left\{1\otimes  x: x\in \ell_2\right\}
 \subset L_2^c(u_c;\ell_2)+L_2^r(u_r;\ell_2) \,.$$
Here $L_2(u_c;\ell_2)=L_2(u_c)\ot_2\el_2$ denotes the
$\el_2$-valued weighted $L_2$-space, and the two weights are
assumed to satisfy
 $$\int_{\O}\min(u_c,\,u_r)d\nu<\8.$$
Probabilistically, we may think of $K_{u_c,\,u_r}$ as obtained by
a family of independent copies of a single variable $v(1)$, where
$v: L_2^c(u_c)+L_2^r(u_r)\to L_1(\N)$ is the complete embedding
using a quasi-free algebra $\N$. Hence the infinite tensor product
$V=v^{\otimes{\infty}}$ yields an embedding of $K_{u_c,\,u_r}$ in
$L_1(\N^{\otimes{\infty}})$. Note that in the classical setting
$\N^{\otimes{\infty}}$ will be either the hyperfinite II$_1$ or
III$_1$ factor. This construction is motivated by the previous
works on $OH$ which admits such a description. We refer to
\cite{ju-OH}, \cite{ju-araki}, \cite{pis-CRp} and \cite{xu-embed}
for more details and related results.

\medskip

The spaces of type $K_{u_c,\,u_r}$ clearly belong to $HQS(C\op
R)$. The first main result of this paper states that almost all
spaces in $HQS(C\op R)$ admit such representations.

\begin{thm}\label{charc}
 Let  $F\in HQS(C\oplus R)$ be infinite dimensional.
Then there exist constants $\l,\g\ge 0$ and two weights $u_c$,
$u_r$ on $\rz$ with respect to Lebesgue measure such that $F$ is
completely isomorphic to the intersection
 $\l C \cap \g R \cap K_{u_c,\,u_r}\,$.
Moreover, under a mild assumption on the following functions
 \beq\label{ff}
 \F_c(n)=
 \big\|\sum_{k=1}^n e_{k 1}\ot e_k\big\|_{C\ot_{\min}F}^2
 \quad\mbox{and}
 \quad \F_r(n)=
 \big\|\sum_{k=1}^n e_{1k}\ot e_k\big\|_{R\ot_{\min}F}^2\,,
 \quad n\in\nz,
 \eeq
where $(e_k)$ denotes an orthonormal basis of $F$, we have
$\l=\g=0$ and the pair $(u_c,\,u_r)$ is uniquely determined, up to
equivalence, by $\F_c$ and $\F_r$.
 \end{thm}

We will call the two functions in \eqref{ff} the \emph{fundamental
functions} (or \emph{sequences}) of $F$ and denote them by $\F_{c,
F}$ and $\F_{r, F}$ if we wish to refer to $F$ explicitly. The
mild assumption mentioned previously is the following: there exist
positive constants $c,\,d$ and $\a,\,\b$ with $0<\a\le \b<1$ such
that
 \beq\label{reg fund}
 c\left(\frac{n}{k}\right)^\a\le
 \frac{\F_c(n)}{\F_c(k)}
 \le d\left(\frac{n}{k}\right)^\b\quad\mbox{and}\quad
 c\left(\frac{n}{k}\right)^\a\le
 \frac{\F_r(n)}{\F_r(k)}
 \le d\left(\frac{n}{k}\right)^\b\,,\quad\forall\; n\ge k\ge1.
 \eeq
In this case $F$ will be called \emph{regular}. We will
systematically extend $\F_c$ and $\F_r$ to continuous functions on
$\rz_+$ (still denoted by the same symbols), for instance,
piecewise linearly. We can even assume, by perturbation, that
$\F_c$ and $\F_r$ are increasing on $\rz_+$. In particular, under
the mild regularity assumption, every homogenous space in
$QS(C\oplus R)$ admits a complete embedding in the predual of the
hyperfinite III$_1$ factor.

\medskip

Theorem \ref{charc} shows that the fundamental functions $\Phi_c$
and $\Phi_r$ completely determine the operator space structure of
a regular $F\in HQS(C\op R)$.  We understand this theorem as a
classification result of nice Hilbertian operator spaces (the
class $HQS(C\oplus R)$) which should play a similar role for
operator spaces as Hilbert spaces do for general Banach spaces.
Although much work is left to be done in this direction, we show
that the fundamental functions $\F_c$ and $\F_r$ do allow us to
calculate fundamental invariants of the operator space $F$.
Indeed, a new feature in operator space theory is the notion of
exactness. An operator space $F$ is called \emph{exact} if
 $$ex(F)=\sup_{E\subset F,\, \dim E<\8}\,
 \inf_{G\subset\mathbb K(\ell_2)}
 d_{cb}(E,\,G)<\8,$$
where $\mathbb K(\el_2)$ denotes the space of compact operators on
$\el_2$. This notion was first investigated by Kirchberg
\cite{kirch-ex}, \cite{kirch-CAR} and \cite{kirch-nonsemi} in the
category of C*-algebras and then by Pisier \cite{pis-ex} for
operator spaces. In fact, a $C^*$-algebra $A$ is exact in the
categorial sense if and only if $ex(A)=1$. Thus knowing that
$ex(F)>1$ implies that no $C^*$-algebra generated by some copies
of $F$ can be exact. For example, the $n$-dimensional operator
Hilbert space $OH_n$ satisfies $ex(OH_n)\sim n^{1/4}$, whereas
$ex(R)=1=ex(C)$.

\begin{thm}\label{ex}
 Let $F\in HQS(C\op R)$ be regular and  $F_n$ an
$n$-dimensional subspace of $F$. Let $\F_c$ and $\F_r$ be the
fundamental functions of $F$. Then
  $$ex(F_n)\sim \left[\frac{n}{\F_{c}(n)}\,
  \F_{r}\big(\frac{\F_{c}(n)}{\F_{r}(n)}\big) +
  \frac{n}{\F_{r}(n)}\,
  \F_{c}\big(\frac{\F_{r}(n)}{\F_{c}(n)}\big)
  \right]^{1/2},$$
where the equivalence constants depend only on the homogeneity
constant of $F$ and the constants in the regularity condition
\eqref{reg fund}.
\end{thm}

Our next result concerns the projection constants of the spaces
$F_n$. Calculating projection constants for classical Banach spaces
has advanced the general knowledge and theory, in particular the
local theory. We think that this particular class of operator spaces
should serve as a testing class for understanding general features
of operator spaces, notably because this class is amenable to
concrete calculations. Recall that the \emph{projection} or
\emph{injectivity constant} of $F$ is defined by
  $$\l_{cb}(F) = \inf\left\{\|P\|_{cb}\,:\, F\subset B(H)
 \textrm{ as subspace, } P:B(H)\to F \textrm{
 projection}\right\}.$$

\begin{thm}\label{inj}
 Under the same assumption of Theorem~\ref{ex} we have
  \be
 &&\frac1{\l_{cb}(F_n)}
 \sim
 \left[\frac1{\F_{c, F}(n)\F_{r, F^*}(n)} +
 \frac1{\F_{r, F}(n)\F_{c, F^*}(n)}\right]^{1/2}+\\
 &&\hskip 1.1cm
 \frac1{\sqrt n}\,\left[\int_1^{\F_{c, F^*}(n)}
 \frac{\F_{r, F}(\F^{-1}_{c, F^*}(t))}
 {\F^{-1}_{c, F^*}(t)}\,dt+
 \int_1^{\F_{r, F^*}(n)}
 \frac{\F_{c, F}(\F^{-1}_{r, F^*}(t))}
 {\F^{-1}_{r, F^*}(t)}\,dt\;+\right.\\
 &&\hskip 2cm\left.\int_1^{\F_{c, F}(n)}
 \frac{\F_{r, F^*}(\F^{-1}_{c, F}(t))}
 {\F^{-1}_{c, F}(t)}\,dt\;+
 \int_1^{\F_{r, F}(n)}
 \frac{\F_{c, F^*}(\F^{-1}_{r, F}(t))}
 {\F^{-1}_{r, F}(t)}\,dt\,\right]^{1/2}.
 \ee
Here $\F^{-1}$ denotes the generalized inverse of a nondecreasing
function $\F$ on $\rz_+$.
 \end{thm}

It is interesting to note that the right hand side above is
symmetric in $F$ and $F^*$. Consequently, for a regular $F\in
HQS(C\op R)$ we have $\l_{cb}(F_n)\sim \l_{cb}(F_n^*)$ (see
Proposition~\ref{inj-sym} for a more general result of this kind).
We will show in Theorem~\ref{dual} below that if $F\in HQS(C\op
R)$ is regular, so is $F^*$ and their fundamental functions are
linked as follows
 $$\F_{c, F}(n)\F_{c, F^*}(n)\sim n \quad\mbox{and}\quad
 \F_{r, F}(n)\F_{r, F^*}(n)\sim n. $$

\medskip

Let $C_p$ denote the column $p$-space, which is the column
subspace of the Schatten $p$-class. It is known that $C_p\in
QS(C\op R)$ (see \cite{ju-OH} and \cite{xu-embed}). It is also
easy to calculate its fundamental functions
 $$\F_{c, C_p}(t)=t^{1/p'}\quad\mbox{and}\quad
 \F_{r, C_p}(t)=t^{1/p}\,,\quad t\in\rz_+\,,$$
where $p'$ denotes the conjugate index of $p$. Thus $C_p$ is
regular for $1<p<\8$, and by the preceding theorems, we
immediately find
 $$ex(C_p^n)\sim n^{1/pp'}$$
and
 $$\l_{cb}(C_p^n)\sim n^{1/\max(p, \;p')}\;
 \textrm{ if } p\neq2
 \quad\textrm{and}\quad
 \l_{cb}(C_2^n)\sim \frac{\sqrt n}{\sqrt{\log(n+1)}}\,.$$
The estimates above for $p\neq 2$ are new. For $p=2$ the estimate
on $ex(C_2^n)$ is due to Pisier (\cite{pis-ex};
see also \cite[Theorem~21.5]{pis-intro}) and that
on $\l_{cb}(C_2^n)$ is the combination of
\cite[Corollary~3.7]{pisshlyak} and \cite[Corollary~4.11]{ju-OH}.

\medskip

One main new feature in our arguments is the use of Orlicz space
techniques. This is the first time that these techniques are
employed in operator space theory. They also allow us to describe
the completely $1$-summing maps in Effros-Ruan's sense between two
spaces $E$ and $F$ in $HQS(C\op R)$.  Let $\Pi_1^o(E,\,F)$ denote
the space of all completely $1$-summing maps from $E$ to $F$,
equipped with the completely $1$-summing norm $\pi_1^o$. Using the
representation Theorem~\ref{charc} we show that $\Pi_1^o(E,\,F)$
isomorphically coincides  with a Schatten-Orlicz class $S_\f$ (see
the beginning of section~\ref{Completely 1-summing maps} for the
definition of $S_\f$). Since $S_\f$ is determined, up to an
equivalent norm, by the fundamental sequence of the underlying
Orlicz function $\f$, this result reduces the determination of
$\Pi_1^o(E,\,F)$ to that of the sequence
($\pi_1^o(\id_n))_{n\ge1}$, where $\id_n$ is the formal identity
from an $n$-dimensional subspace of $E$ to another one of $F$.
Choose two  orthonormal bases $(e_k)$ and $(f_k)$ of $E$ and $F$,
respectively. Then $\id_n$ is the map such that $\id_n(e_k)=f_k$
for $k\le n$ and $\id_n(e_k)=0$ for $k>n$. The homogeneity of $E$
and $F$ shows that $\pi_1^o(\id_n)$ is independent of particular
choice of $(e_k)$ and $(f_k)$.

Only with the help of this reduction result we can determine the
whole space $\Pi_1^o(E,\,F)$. In many situations it is hard or
even impossible to estimate $\pi_1^o(u)$ for any $u:E\to F$ but it
is relatively easier to determine $\pi_1^o(\id_n)$ via concrete
integral calculations. This is indeed the case for $E=F=OH$. Then
$\pi_1^o(\id_n\,:\, OH_n\to OH_n)$ was determined  by the first
named author in \cite{ju-OH}: $\pi_1^o(\id_n\,:\, OH_n\to
OH_n)\sim\sqrt{n\log(n +1)}$ uniformly in $n\in\nz$. Note that
this estimate is equivalent to that of $\g_{cb}(OH_n)$ in
Theorem~\ref{inj}. Thus the previous reduction result implies that
$\Pi_1^o(OH)=S_\psi$, where $S_\psi$ is the Schatten-Orlicz class
associated to the Orlicz function $\psi$ defined by
$\psi(t)=t^2\log(t+1/t)$.  This result is in strong contrast with
the corresponding result in Banach space theory. Recall that a map
$u: X\to Y$ between two Banach spaces is $1$-summing if there
exists a constant $\l>0$ such that for all finite sequences
$(x_k)\subset X$
 $$\sum_k\|u(x_k)\|\le \l\,\sup\big\{\sum_k|\xi(x_k)|\;:\; \|\xi\|\le
 1,\; \xi\in X^*\big\}.$$
It is  well known that a map $u$ on a Hilbert space $H$ is
$1$-summing if and only if $u$ is a Hilbert-Schmidt operator (see,
for instance, \cite{pis-fact}).

\medskip

If both $E$ and $F$ are regular, we can do much better. In this
case, we have an explicit formula for $\pi_1^o(\id_n)$ in terms of
the fundamental functions of $E$ and $F$.

\begin{thm}\label{c1s id explicit}
 Let $E, F\in HQS(C\op R)$ be regular. Then
$\Pi_1^o(E,\,F)=S_\f$ for some Orlicz function $\f$ whose
fundamental sequence $(\f_n)_{n\ge1}$ is given by
 \be
 &&\f_n^2\sim
 \F_{c, E^*}(n)\F_{r, F}(n)\;+\;
 \F_{r, E^*}(n)\F_{c, F}(n)\;+\\
 &&\hskip .8cm
 n\Big[\int_1^{\F_{c, E^*}(n)}\frac{\F_{r, F}(\F^{-1}_{c, E^*}(t))}
 {\F^{-1}_{c, E^*}(t)}\,dt+
 \int_1^{\F_{r, E^*}(n)}
 \frac{\F_{c, F}(\F^{-1}_{r, E^*}(t))}
 {\F^{-1}_{r, E^*}(t)}\,dt+\\
 &&\hskip 1.2cm \int_1^{\F_{c, F}(n)}
 \frac{\F_{r, E^*}(\F^{-1}_{c, F}(t))}
 {\F^{-1}_{c, F}(t)}\,dt\;+
 \int_1^{\F_{r, F}(n)}
 \frac{\F_{c, E^*}(\F^{-1}_{r, F}(t))}
 {\F^{-1}_{r, F}(t)}\,dt\,\Big],
 \ee
where the equivalence constants depend only on the homogeneity and
regularity constants of $E$ and $F$.
 \end{thm}

Applying this theorem to the column $p$-spaces, we immediately
determine the whole space $\Pi_1^o(C_p,\, C_q)$. In particular, we
recover the result on $\pi_1^o(\id_n\,:\, C_2^n\to C_2^n)$ of
\cite{ju-OH} quoted above  by simpler arguments which avoid some
tedious integral calculations (for instance, we do not use the
duality argument of \cite{ju-OH}).

\medskip

The paper is organized as follows. After a preliminary section on
operator space theory, we prove in section~\ref{Representations of
homogeneous spaces} the first part of Theorem~\ref{charc}, which
is reformulated as Theorem~\ref{hom rep}. The uniqueness part of
Theorem~\ref{charc} is proved in section~\ref{Homogeneous spaces
satisfying a regularity condition}. There the Orlicz space
techniques mentioned previously appear for the first time (see the
proof of Lemma~\ref{orlicz cr}).  Section~\ref{Completely
1-summing maps} deals with completely $1$-summing maps. The first
main result there states that for $E$ and $F$ in $HQS(C\op R)$ we
have $\Pi_1^o(E,\, F)=S_\f$ isomorphically for some Orlicz
function $\f$. Another main result is Theorem~\ref{c1s id
explicit}. Section~\ref{Injectivity and exactness} concerns the
projection constants and exactness constants. There
Theorems~\ref{inj} and \ref{ex} are proved. In the last section we
apply all these results to the column spaces $C_p$ and their sums
and intersections with the row spaces $R_p$. Consequently, we
determine all previous objects for these spaces.

\medskip

The techniques developed in this paper allow to obtain natural
operator space structures on certain Schatten-Orlicz spaces. They
also permit to deal with completely $p$-summing maps in Pisier's
sense. These subjects will be pursued elsewhere. We refer to Yew's
paper \cite{yew} for the study of completely $p$-summing maps on
$OH$.

All spaces considered in this paper will be separable and infinite
dimensional, unless explicitly stated otherwise. The letter $c$
will often denote a universal constant. We will frequently use the
notation $A\sim_cB$ to abbreviate the two-sided inequality
$c^{-1}B\le A\le cB$.


\section{Preliminaries}
 \label{Preliminaries}


In this section we collect some preliminaries necessary to the
whole paper. We will use standard notions and notation from
operator space theory. Our references are \cite{er-book} and
\cite{pis-intro}. An operator space $E$ is called {\it
homogeneous} if there exists a constant $\l$ such that every
bounded map $u$ on $E$ is completely bounded  and
$\|u\|_{cb}\le\l\|u\|$. In this case we also say that $E$ is
$\l$-homogeneous. $E$ is called {\it Hilbertian} if $E$ is
isomorphic to a Hilbert space (at the Banach space level). If we
wish to emphasize the isomorphism constant $\l$ between $E$ and
the Hilbert space, we say that $E$ is $\l$-Hilbertian.

The Schatten classes $S_p$ will be frequently used in this paper.
Recall that $S_1$ is the trace class, $S_2$ the Hilbert-Schmidt
class and $S_\8=B(\el_2)$.  These spaces are equipped with their
natural operator space strictures as introduced in \cite{pis-ast}.
We will also need their vector-valued versions. Let $E$ be an
operator space. Define $S_\8[E]=S_\8\ot_{\min}E$ and
$S_1[E]=S_1\wh\ot E$. Here $\ot_{\min}$ and $\wh\ot$ denote,
respectively,  the minimal (injective) and projective tensor
products in the category of operator spaces. For $1<p<\8$ the
space $S_p[E]$ is defined as the complex interpolation space
$(S_\8[E],\; S_1[E])_{1/p}$. Only $S_\8[E]$ and $S_1[E]$ will be
needed later. We refer to \cite{pis-ast} for more information.

We now recall the direct sum of two operator spaces $E$ and $F$.
Let $1\le p\le\8$. $E\op_p F$ denotes the direct sum of $E$ and
$F$ in the $\el_p$-sense (see \cite{pis-ast}). Note that $E\op_pF$
is completely isomorphic to $E\op_qF$ for any $1\le q\le\8$ with
universal constants. This allows us to drop the subscript $p$ from
$E\op_pF$, a convention adopted throughout the paper. Note that if
both $E$ and $F$ are Hilbertian, so is $E\op F$. However, in the
$1$-Hilbertian case, only $E\op_2F$ is $1$-Hilbertian.

We will need to consider the sum and intersection of a compatible
couple $(E, \, F)$ of operator spaces. The intersection $E\cap F$
is the diagonal subspace of $E\op F$ and the sum $E+F$ is
 $$E+F=\{a+b\,:\, a\in E,\,b\in F\}.$$
Note that $E+F$ is the quotient of $E\op F$ by the subspace
$\{(a,b)\,:\, a+b=0\}$.

The column and row spaces, $C$ and $R$, are the two major objects
of this paper. Recall that $C$ and $R$ are respectively the
(first) column and row subspaces of $S_\8$. More precisely, $C$ is
the subspace of $S_\8$ consisting of matrices whose all entries
but those in the first column  vanish. As Banach spaces both $C$
and $R$ are isometric to $\el_2$, so they are $1$-Hilbertian. This
allows us to identify both $C$ and $R$ with $\el_2$ at the Banach
space level. Accordingly, we will often identify the canonical
bases $(e_{k1})$ of $C$ and $(e_{1k})$ of $R$ with $(e_k)$ of
$\ell_2$. It is easy to see that the operator space structures of
$C$ and $R$ are determined as follows. Let $(x_k)$ be a finite
sequence in $S_\8$. Then
 $$\big\|\sum_kx_k\ot e_k\big\|_{S_\8[C]}
 =\big\|\sum_k x_k^*x_k\big\|_\8^{1/2}\,,\quad
 \big\|\sum_kx_k\ot e_k\big\|_{S_\8[R]}
 =\big\|\sum_k x_kx_k^*\big\|_\8^{1/2}\,.$$
This implies that $C$ and $R$ are $1$-homogeneous. More generally,
if $H$ is a Hilbert space, the column and row spaces based on $H$
are $H^c=B(\cz,\, H)$ and $H^r=B(H,\, \cz)$, respectively. If $H$
is separable and infinite dimensional, we recover $C$ and $R$. On
the other hand, if $\dim H=n<\8$, we get $C^n$ and $R^n$, the
$n$-dimensional versions of $C$ and $R$.

We will be interested only in the homogeneous spaces in $QS(C\op
R)$. Here given an operator space $E$ we use $QS(E)$ to denote the
family of all quotients of subspaces of $E$. This coincides with
the family of all subspaces of quotients of $E$. The subfamily of
homogeneous spaces of $QS(C\op R)$ is denoted by $HQS(C\op R)$. We
will study several properties of these spaces.  It is well-known
that Pisier's operator Hilbert space $OH$ constructed in
\cite{pis-oh} belongs to $HQS(C\op R)$ (see \cite[Exercice
7.9]{pis-intro} and also \cite{pis-CRp}). More generally, the
column $p$-space $C_p$ belongs to $HQS(C\op R)$ too (see
\cite{ju-OH} and \cite{xu-embed}).

Recall that $C_p$ and $R_p$ denote the (first) column and row
subspaces of $S_p$. Their $n$-dimensional versions are denoted by
$C_p^n$ and $R_p^n$, respectively. Note that $C_\8$ and $R_\8$ are
just $C$ and $R$. On the other hand, $C_2$ and $R_2$ coincide
completely isometrically with $OH$. Like $C$ and $R$, $C_p$ and
$R_p$ are also 1-homogenous and 1-Hilbertian.  We will also
identify $C_p$ and $R_p$ with $\el_2$ as Banach spaces and use
$(e_k)$ to denote their common canonical basis. We have the
following completely isometric identities: for any $1\le p\le\8$
 \be
 (C_p)^*\cong C_{p'}\cong R_p\quad\mbox{and}\quad
 (R_p)^*\cong R_{p'}\cong C_p,
 \ee
where $p'$ denotes the index conjugate to $p$. $C_p$ and $R_p$ can
be also defined via interpolation from $C$ and $R$. We view $(C,\,
R)$ as a compatible couple by identifying both of them with
$\ell_2$ (at the Banach space level). Then
 $$C_p=(C,\; R)_{1/p}=(C_\8,\; C_1)_{1/p}\quad\mbox{and}\quad
 R_p=(R,\; C)_{1/p}=(R_\8,\; R_1)_{1/p}\ .$$
We refer to \cite{pis-oh} and \cite{pis-ast} for all these
elementary facts.

If $E$ is an operator space, $C_p[E]$ (resp. $R_p[E]$) denotes the
closure of $C_p\ot E$ (resp. $R_p\ot E$) in $S_p[E]$. Thus
$C[E]=C\ot_{\min} E$ and $C_1[E]=C_1\wh\ot E$.

We end this preliminary section by introducing the notion of
completely $1$-summing maps. Let $x: E\to F$ be a  map between two
operator spaces. $x$ is called {\it completely $1$-summing} if the
map $\id\ot x$ is bounded from $S_1\ot_{\min} E$ to $S_1[F]$. In
this case we define $\pi_1^o(x)$ to be the norm of  $\id\ot x$ and
call it the {\it completely $1$-summing norm of $x$}. The space of
all completely $1$-summing maps from $E$ to $F$ is denoted by
$\Pi_1^o(E,\, F)$ and equipped with the completely $1$-summing
norm. This is a Banach space. It is easy to check that
$\Pi_1^o(E,\, F)$ is an ideal in the following sense. Let $E_1$
and $F_1$ be two other operator spaces and let $y\in CB(E_1,\, E)$
and $z\in CB(F,\, F_1)$. Then $zxy\in\Pi_1^o(E_1,\, F_1)$ for any
$x\in\Pi_1^o(E,\, F)$ and
$\pi_1^o(zxy)\le\|z\|_{cb}\,\pi_1^o(x)\,\|y\|_{cb}$. We refer to
\cite{er-gro-pist} and \cite{pis-ast} for more information.


\section{Representations  of homogeneous spaces in $QS(C\op R)$}
 \label{Representations  of homogeneous spaces}


In this section we consider spaces in $HQS(C\op R)$. The main
result is a representation theorem for these spaces (the first
part of Theorem~\ref{charc}). Let $(\O, \nu)$ be a measure space
and $u$ a weight on $\O$ (i.e., a nonnegative measurable
function). We denote by $L_2(\O, u)$ the corresponding weighted
$L_2$-space whose norm is given by
 $$\|f\|_{L_2(\O, u)}=\big(\int_\O |f|^2u d\nu\big)^{1/2}\,.$$
Similarly, given a Banach space $X$ we define the space $L_2(\O,u;
X)$ of functions on $\O$ with values in $X$; the norm of
$L_2(\O,u;X)$ is defined as above by replacing the absolute value
by the norm of $X$. We will need, however, only the case
$X=\el_2$. Then $L_2(\O,u; \el_2)$ is again a Hilbert space.
$L_2(\O, u)$ and $L_2(\O,u;X)$ will be denoted simply by $L_2(u)$
and $L_2(u;X)$, respectively, if no confusion can occur. We will
denote by $L_2^c(u;\el_2)$ (resp. $L_2^r(u;\el_2)$)  the column
(resp. row) operator space based on $L_2(u;\el_2)$.

Let $(u_c,\,u_r)$ be a pair of weights on $\O$ such that
 \beq\label{w-con}
 \int_\O\min(u_c,\; u_r)\,d\nu<\8.
 \eeq
We will call this the {\it weight condition} and will always
assume it whenever a pair of weights is considered.   We view
$(L_2^c(u_c;\el_2)\,,\, L_2^r(u_r;\el_2))$ as a compatible pair in
the natural way. Then define
 \beq\label{def of G}
 G_{u_c,\,u_r}=L_2^c(u_c;\el_2)+ L_2^r(u_r;\el_2)\,.
 \eeq
This is the quotient of $L_2^c(u_c;\el_2)\op L_2^r(u_r;\el_2)$ by
the subspace of all $\el_2$-valued functions $(a, b)$ such that
$a+b=0$ a.e. on $\O$; so $G_{u_c,\,u_r}$ can be viewed as a
quotient of $C\op R$. Let $K_{u_c,\,u_r}$ be the subspace of
constant functions of $G_{u_c,\,u_r}$. Identifying constant
functions with vectors of $\el_2$, we easily check that
$K_{u_c,\,u_r}$ coincides (isomorphically) with $\el_2$ as Banach
spaces. Moreover, the isomorphism becomes an isometry if the
underlying sum is in the $\el_2$-sense. Indeed, given $x=(x_k)\in
K_{u_c,\,u_r}$ we have
 \be
 \|x\|_{K_{u_c,\,u_r}}^2
 &=&\inf\big\{\sum_k\int_\O \big(|a_k|^2\,u_c d\nu+
 |b_k|^2\,u_r \big)d\nu\,:\, x_k=a_k(\o)+b_k(\o) \;
 \textrm{a.e.}\big\}\\
 &=&\sum_{k}\int_\O\inf_{0\le t\le 1}(t^2u_c
 +(1-t)^2u_r )\,|x_k|^2\,d\nu
 =\int_\O\frac{u_cu_r}{u_c +u_r }\,\sum_k|x_k|^2\,.
 \ee
The last integral  is finite thanks to \eqref{w-con}. In the
sequel we will often identify $K_{u_c,\,u_r}$ and $\el_2$ at the
Banach level and use $(e_k)$ to denote the canonical basis of
$K_{u_c,\,u_r}$ too.

The operator space structure of $K_{u_c,\,u_r}$ is given as
follows. For any finite sequence $(x_k)\subset S_\8$
 \be
 \big\|\sum_kx_k\ot e_k\big\|_{S_\8[K_{u_c,\,u_r}]}=\inf\Big\{
 \big\|\sum_k\int_\O a_k^*a_k\,u_c d\nu\big\|_{S_\8}^{1/2}+
 \big\|\sum_k\int_\O b_kb_k^*\,u_r d\nu\big\|_{S_\8}^{1/2}
 \Big\},
 \ee
where the infimum runs over all decompositions
$x_k=a_k(\o)+b_k(\o)$ a.e. on $\O$ with $a_k\in S_\8[L_2^c(u_c)]$
and $b_k\in S_\8[L_2^r(\mu)]$. We will show that these spaces
$K_{u_c,\,u_r}$ are the nontrivial building blocks of spaces in
$HQS(C\op R)$.

\begin{rk}\label{positive weight}
 If $u_c$ or $u_r$ vanishes  on a subset $A\subset\O$,
then $A$ does not contribute to $K_{u_c,\,u_r}$. Namely, the space
$K_{u_c,\,u_r}$ constructed over $\O$ is the same as that over
$\O\setminus A$. Thus all weights in the sequel will be assumed
strictly positive unless explicitly stated otherwise.
 \end{rk}

It is sometimes convenient to work with the discrete analogue of
$K_{u_c,\,u_r}$, i.e., when $(\O, \nu)$ is a discrete measure
space. Consider, for instance, the case where $\O=\nz$ and $\nu$
is the counting measure. Then the two weights $u_c$ and $u_r$
become two positive sequences $(u_c(j))_{j\ge1}$ and
$(u_r(j))_{j\ge1}$ satisfying the following weight condition
 \beq\label{w-con dis}
 \sum_j\min(u_c(j),\; u_r(j))<\8.
 \eeq
The space $G_{u_c,\,u_r}$ is now given by
 $$G_{u_c,\,u_r}=\el_2^c(u_c;\el_2) + \el_2^r(u_r;\el_2)\,.$$
$K_{u_c,\,u_r}$ is the subspace of $G_{u_c,\,u_r}$ consisting of
constant sequences.

\medskip

By standard arguments it is easy to transfer the continuous case
to the discrete one and vice versa. More precisely, we have the
following

\begin{prop}\label{cont-dis}
 \begin{enumerate}[\rm(i)]
 \item Let $(u_c,\,u_r)$ be a pair of weights on a measure space
$(\O,\nu)$ verifying \eqref{w-con}. Let $(d_j)$ be a positive
sequence and $(A_j)$ a partition of $\O$ such that
 $$d_j\le\, \frac{u_r}{u_c}\,\le c\, d_j\quad\mbox{on}\quad A_j$$
for some positive constant $c$. Let $\wt u_c=(\wt u_c(j))$ and
$\wt u_r=(\wt u_r(j))$ be two sequences defined by
 $$\wt u_c(j)=\int_{A_j}u_c\,d\nu
 \quad\mbox{and}\quad \wt u_r(j)=d_j\,\wt u_c(j).$$
Then $(\wt u_c,\,\wt u_r)$ satisfies \eqref{w-con dis} and
$K_{u_c,\,u_r}$ is completely isomorphic to $K_{\wt u_c,\,\wt
u_r}$.
 \item Conversely, given two positive sequences $\wt u_c$
and $\wt u_r$ verifying \eqref{w-con dis} define two weights $u_c$
and $u_r$ on $\rz_+$ by
 {\rm $$u_c=\sum_{j\ge 1} \wt u_c(j)\, \un_{(j-1,\; j]}
 \quad\mbox{\em and}\quad
 u_r=\sum_{j\ge 1} \wt u_r(j)\, \un_{(j-1,\; j]}\,.$$}
Then $u_c$ and $u_r$ satisfy \eqref{w-con} and $K_{\wt u_c,\,\wt
u_r}$ is completely isomorphic to $K_{u_c,\,u_r}$, where $\rz_+$
is equipped with Lebesgue measure.
 \end{enumerate}
 \end{prop}

\pf (i) We have
 \be
 \sum_j\min(\wt u_c(j)\,,\; \wt u_r(j))
 &=&\sum_j\int_{A_j}\min(1,\; d_j) u_c d\,\nu\\
 &\le& \sum_j\int_{A_j}\min(u_c,\; u_r)d\,\nu
 =\int_\O\min(u_c,\; u_r) d\,\nu<\8.
 \ee
Define a map $T$ by
 $$T(f)=\big(\frac{1}{\wt u_c(j)}\,
 \int_{A_j} f u_c\,d\nu\big)_{j\ge1}\,.$$
Then it is  easy to check that $T$ is a contraction from
$L_2(u_c)$ to $\el_2(\wt u_c)$  as well as from $L_2(u_r)$ to
$\el_2(\wt u_r)$. It follows, by homogeneity,  that $T\ot\id$ is
completely contractive from $L_2^c(u_c;\el_2)$ to $\el_2^c(\wt
u_c;\el_2)$ and from $L_2^r(u_r;\el_2)$ to $\el_2^r(\wt
u_r;\el_2)$. Now let $(x_k)$ be a finite sequence in $S_\8$ and
consider a decomposition $x_k=a_k(\o)+b_k(\o)$ a.e. on $\O$. Then
$x_k=\wt a_{jk}+\wt b_{jk}$ for all $j$, where $\wt a_k =T(a_k)$
and $\wt b_k =T(b_k)$. Moreover,
 $$\big\|\sum_{j,k} \wt u_c(j)\,\wt a_{jk}^*\wt a_{jk}\big\|
 \le \big\|\sum_k\int_\O a_k^*a_ku_cd\nu\big\|\,,\quad
 \big\|\sum_{j,k}\wt u_r(j)\, \wt b_{jk}\wt b_{jk}^*\big\|
 \le \big\|\sum_k\int_\O b_kb_k^*u_rd\nu\big\|\,.$$
We thus deduce
 $$\big\|\sum_kx_k\ot e_k\big\|_{S_\8[K_{\wt u_c,\,\wt u_r}]}\le
 \big\|\sum_kx_k\ot e_k\big\|_{S_\8[K_{u_c,\,u_r}]}\,.$$
The converse inequality is proved similarly by using the map $T'$
defined by
 $$T'(x)=\sum_jx_j\un_{A_j}\,.$$
Indeed, $T'$ is contractive from $\el_2(\wt u_c)$ to $L_2(u_c)$
and bounded from $\el_2(\wt u_r)$ to $L_2(u_r)$ with norm
$\le\sqrt c$. Then as above, we deduce the missing converse
inequality. Therefore, $K_{u_c,\,u_r}$ is completely isomorphic to
$K_{\wt u_c,\,\wt u_r}$. This shows (i). The proof of (ii)  is
similar and thus omitted.\cqd

\medskip

The following theorem shows that except $C$, $R$ and $C\cap R$ all
homogeneous spaces in $QS(C\op R)$ are of the form $K_{u_c,\,u_r}$
for some sequences $u_c$ and $u_r$. Recall that if $E$ is an
operator space and $\l$ a positive constant,  $\l E$ denotes the
same operator space as $E$ but with norm equal to $\l$ times that
of $E$. For convenience we also set $\l E=\{0\}$ if $\l=0$. In the
latter case the intersection $\l C\cap\, \g R\,\cap K_{u_c,\,u_r}$
below is understood as $\g R\,\cap K_{u_c,\,u_r}$.

\begin{thm}\label{hom rep}
 Let $F$ be an infinite dimensional space in
$HQS(C\op R)$. Then there exist two constants $\l, \g\in[0,\;1]$
and two positive sequences $u_c=(u_c(j))_{j\ge1},\;
u_r=(u_r(j))_{j\ge1}$ such that $u_c$ and $u_r$ satisfy
\eqref{w-con dis} and such that $F$ is completely isomorphic to
$\l C\cap\, \g R\,\cap K_{u_c,\,u_r}$. Moreover, the relevant
constant depends only on the homogeneity constant of $F$.
 \end{thm}

We start the proof of the theorem by some preparations.  Let
$S\subset C\op R$ be a closed subspace such that $F$ is a subspace
of the quotient $(C\op R)/S$. By the decomposition theorem of
\cite{xu-embed}, we find four subspaces $H_j\subset \el_2$, $0\le
j\le 3$ and an injective closed densely defined operator $T:
H_2\to H_3$ of dense range such that
 $$S= H_0^c\op H_1^r\op \Ga(T),$$
where $\Ga(T)=\{(x,\, Tx)\,:\, x\in{\rm Dom}(T)\}$ is the graph of
$T$, viewed as a subspace of $C\op R$. Assume that all direct sums
are in the $\el_2$-sense. Then the previous decomposition of $S$
is orthogonal. This implies $H_0\perp H_2$ and $H_1\perp H_3$. On
the other hand, writing the polar decomposition of $T$ and using
homogeneity, we can assume that $H_2=H_3$ and $T$ is a positive
operator on $H_2$. Next, using the spectral decomposition of $T$
and by approximation, we can further assume that $T$ has only
point spectrum, i.e., $H_2$ has an orthonormal basis consisting of
eigenvectors of $T$. Finally, by homogeneity once more, we may fix
an orthonormal basis in each $H_j$ which consists of vectors in
the canonical basis $(e_k)$ of $\el_2$. Moreover, the basis of
$H_2$ is formed of eigenvectors of $T$. This choice of $S$ will be
fixed in the sequel. It is to ensure the invariance of $S$ by any
diagonal operator on $\el_2$.

Fix an orthonormal basis $(f_k)$ of $F$. Then $(f_k)$ is
\emph{completely symmetric} in the following sense: there exists a
constant $\l$ (majorized by the homogeneity constant of $F$) such
that
 $$ \|\sum_k \e_k a_{\pi(k)} \ot f_k \| \le\l\,
  \|\sum_k  a_{k } \ot f_k \| $$
holds for all finite sequences $(a_k)\subset S_\8$, $\e_k=\pm 1$
and permutations $\pi$ on $\nz$. This is equivalent to say that
$(f_k)$ is completely equivalent to $(f_{\pi(k)})$ for any
permutation $\pi$ of $\nz$. By the complete equivalence of two
bases $(f_k)$ of $F$ and $(g_k)$ of $G$ we mean that the map
$f_k\mapsto g_k$ extends to a complete isomorphism from $F$ onto
$G$. We will assume, for simplicity, that $F$ is $1$-homogeneous.
Set $C(\nz^2)=\el_2(\nz^2)^c$ and $R(\nz^2)=\el_2(\nz^2)^r$.

\begin{lem}\label{disj support}
 There exists a subsequence $(f_{n_k})$ of $(f_k)$ such that
$(f_k)$ is completely equivalent to the basic sequence
$(f_{n_k}\ot e_k)$ in $(C(\nz^2)\op R(\nz^2))/\el_2(S)$.
 \end{lem}

\pf Let $q: C\op R\to (C\op R)/S$ be the quotient map and $(\d_k)$
a positive sequence. Since $(f_k)$ is a basic sequence in $(C\op
R)/S$ and weakly converges  to $0$, by a standard perturbation
argument we find two increasing sequences $(n_k)$ and $(m_k)$ of
positive integers such that
 $$\big\|f_{n_k} - \sum_{m_k\le i<m_{k+1}}
 \big(\a_{i}q(e_{i1})+\b_{i}q(e_{1i})\big) \big\|<\d_k$$
for some  $a_{i}, \b_i\in\cz$. It follows that if the $\d_k$ are
sufficiently small,  $(f_{n_k})$ is $2$-completely equivalent to
$(\wt f_{n_k})$, where
 $$\wt f_{n_k}= \sum_{m_k\le i<m_{k+1}}
 \big(\a_{i}q(e_{i1})+\b_{i}q(e_{1i})\big).$$
Replacing  $(f_{n_k})$ by $(\wt f_{n_k})$ if necessary, we may
assume $f_{n_k}=\wt f_{n_k}$ for all $k$. In this case we will say
that the $f_{n_k}$'s have disjoint supports. On the other hand,
the complete symmetry of $(f_k)$ implies that $(f_k)$ is
completely equivalent to $(f_{n_k})$. Thus it remains to show that
$(f_{n_k})$ is completely equivalent to $(f_{n_k}\ot e_k)$.

Given a finite sequence $(x_k)\subset S_\8$ set
 $$x=\sum_kx_k\ot f_{n_k}\quad\mbox{and}\quad
 \wt x=\sum_kx_k\ot f_{n_k}\ot e_k\,.$$
We are going to show
 \beq\label{xx}
 \|x\|_{S_\8[(C\op R)/S]}=
 \|\wt x\|_{S_\8[(C(\nz^2)\op R(\nz^2))/\el_2(S)]}\,.
 \eeq
Let $a\in S_\8[C]$ and $b\in S_\8[R]$ such that
 $x=\id\ot q(a,\, b)$.
Now let $P_k$ be the projection onto the interval $[m_k,\;
m_{k+1})$. Namely, $P_k$ is the projection on $\el_2$  such that
$P_k(e_i)=e_i$ if $m_k\le i<m_{k+1}$ and $P_k(e_i)=0$ otherwise.
$P_k$ is viewed as projections on both $C$ and $R$ (recalling that
we identify $(e_{i1})$ and $(e_{1i})$ with $(e_i)$). Since
$P_k(S)\subset S$, thanks to the invariance of $S$ by diagonal
operators, $P_k$ induces a projection $\wt P_k$ on $(C\op R)/S$.
Then $qP_k=\wt P_kq$. Thus we find
 $$x_k\ot f_{n_k}=\id\ot \wt P_k(x)
 =\id\ot q\big(\id\ot P_k(a),\,\id\ot P_k(b)\big)
 \;{\mathop=^{\rm def}}\;\id\ot q(\wt a_k,\,
 \wt b_k).$$
Therefore
 $$\wt x=\id\ot q\ot\id( \sum_k\wt a_k\ot e_k,\; \sum_k\wt b_k\ot e_k)
  \;{\mathop=^{\rm def}}\;\id\ot q\ot\id(\wt a,\, \wt b)\,.$$
It is clear that
 $$\|\wt a\|_{S_\8[C(\nz^2)]}\le \|a\|_{S_\8[C]}
 \quad\mbox{and}\quad
 \|\wt b\|_{S_\8[R(\nz^2)]}\le \|b\|_{S_\8[R]}\,;$$
whence
 $$\|\wt x\|_{S_\8[(C(\nz^2)\op R(\nz^2))/\el_2(S)]}\le
 \|x\|_{S_\8[(C\op R)/S]}\,.$$
Conversely,  let $\wt a\in S_\8[C(\nz^2)]$, $\wt b\in
S_\8[R(\nz^2)]$ such that $\wt x=\id\ot q\ot\id(\wt a, \wt b)$. As
before, let $Q_k$ be the projection onto the subset $[m_k,\,
m_{k+1})\times\{k\}$ (considered as a projection on
$\el_2(\nz^2)$). Put $a_k=Q_k(\wt a)$ and $b_k=Q_k(\wt b)$. We
regard $a_k$ and $b_k$ as elements in $S_\8[C]$ and $S_\8[R]$,
respectively. Then
 $$x_k\ot f_{n_k}=\id\ot q(a_k, b_k).$$
Thus
 $$x=\sum_kx_k\ot f_{n_k}=
 \id\ot q\big(\sum_ka_k,\;  \sum_kb_k\big)
 \;{\mathop=^{\rm def}}\;\id\ot q(a, b).$$
It is easy to see
 $$\|a\|_{S_\8[C]}\le\|\wt a\|_{S_\8[C(\nz^2)]}
 \quad\mbox{and}\quad
  \|b\|_{S_\8[R]}\le\|\wt b\|_{S_\8[R(\nz^2)]}\,;$$
so
 $$\|x\|_{S_\8[(C\op R)/S]}\le
 \|\wt x\|_{S_\8[(C(\nz^2)\op R(\nz^2))/\el_2(S)]}\,.$$
Thus \eqref{xx} is proved. Therefore, $(f_{n_k})$ is
$1$-completely equivalent to $(f_{n_k}\ot e_k)$. \cqd

\begin{rk}\label{disjoint}
 The preceding proof shows that if the $f_k$'s have disjoint supports,
then $(f_{k})$ is $1$-completely equivalent to $(f_{k}\ot e_k)$.
 \end{rk}

We will use ultraproducts of operator spaces (see
\cite{pis-intro}). If $\U$ is a free ultrafilter on $\nz$ and $E$
an operator space, we denote by $E^\U$ the ultrapower $\prod_\U
E$. Recall that $E$ is naturally identified as a subspace of
$E^\U$. Let $H$ be the Hilbert space ultrapower of $\el_2$ (in the
category of Banach spaces). Then
 $$\big(\frac{C\op R}{S}\big)^\U\,=\,\frac{H^c\op
 H^r}{{S^\U}}\,.$$
Consequently, $F^\U$ is a subspace of $(H^c\op H^r)/S^\U$. Let
$\wt f$ be the element of $F^\U$ represented by the basis $(f_k)$.

\begin{lem}\label{ultrapower}
 For any finite sequence $(x_k)\subset S_\8$
 we have
 $$\big\|\sum_{k} x_k\ot f_k\big\|_{S_\8[F]}\sim_c
 \big\|\sum_{k} x_k\ot \wt f\ot e_k
 \big\|_{S_\8[((\el_2(H))^c\op (\el_2(H))^r)/\el_2(S^\U)]}\,.$$
 \end{lem}

\pf By Lemma \ref{disj support}, we can assume that the supports
of the $f_k$ are disjoint. Then by the complete symmetry of
$(f_k)$ and Lemma \ref{disj support} (and also Remark
\ref{disjoint}), we find
 $$\big\|\sum_{k=1}^m x_k\ot f_k\big\|_{S_\8[F]}=
 \big\|\sum_{k=1}^m x_k\ot f_{n_k}\ot e_k
 \big\|_{S_\8[(C(\nz^2)\op R(\nz^2))/\el_2(S)]}$$
for any $x_1,..., x_m\in S_\8$ and any distinct positive integers
$n_1, ..., n_m$. We can assume that the $x_k$ are finite matrices
and that all direct sums are taken in the $\el_2$-sense. Then we
find
 $$\big\|\sum_{k=1}^m x_k\ot f_k\big\|_{S_\8[F]}=
 \lim_{n_1, \,\U}\cdots \lim_{n_m, \,\U}
 \big\|\sum_{k=1}^m x_k\ot f_{n_k}\ot e_k
 \big\|_{S_\8[(C(\nz^2)\op R(\nz^2))/\el_2(S)]}\,.$$
This implies the desired assertion. \cqd

\medskip

\n{\it Proof of Theorem~\ref{hom rep}.} Let $(S^\U)^\perp$ be the
orthogonal complement of $S^\U$ in $(H^c\op H^r)^*=\bar H^r\op
\bar H^c$. Then
 $$\big(\frac{H^c\op H^r}{S^\U}\big)^*
 =(S^\U)^\perp\,.$$
On the other hand, by \cite{xu-embed}  we find four subspaces
$\bar K_j\subset \bar H$, $0\le j\le 3$ and an injective closed
densely defined operator $D: \bar K_2\to \bar K_3$ of dense range
such that
 $$(S^\U)^\perp= \bar K_0^r\op \bar K_1^c\op \Ga(D).$$
By the discussion following Theorem \ref{hom rep}, we can assume
that $K_2=K_3=K$, $D$ is a positive operator on $\bar K$ and $\bar
K$ has an orthonormal basis $(\bar g_i)_{i\in I}$ of eigenvectors
of $D$. The eigenvalue associated to $\bar g_i$ is denoted by
$d_i$. We then deduce that
 \be
 \frac{H^c\op H^r}{S^\U}
 &=&\big((S^\U)^\perp\big)^*
 =(\bar K_0^r\op \bar K_1^c)^*\op \Ga(D)^*\\
 &=&K_0^c\op K_1^r
 \op \frac{K^c\op K^r}{\Ga(D)^\perp}\,,
 \ee
where $\Ga(D)^\perp=\{(-Dy,\; y)\,:\, y\in{\rm Dom}(D)\}$. Note
that $\Ga(D)^\perp$ is the closed linear span of all $(-d_ig_i,\;
g_i)$ in $K\op K$. Thus the operator space structure of $(K^c\op
K^r)/\Ga(D)^\perp$ is determined as follows. Let $q: K^c\op K^r\to
(K^c\op K^r)/\Ga(D)^\perp$ be the quotient map. Then $(q(g_i,
0))_{i\in I}$ is a basis of $(K^c\op K^r)/\Ga(D)^\perp$. For any
finite family $(x_i)\subset S_\8$
 $$\big\|\sum_i x_i\ot q(g_i,0)\big\|_{S_\8[(K^c\op K^r)/\Ga(D)^\perp]}
 =\inf_{x_i=a_i+b_i}\max\big(\big\|\sum_ia_i^*a_i\big\|^{1/2},\;
 \big\|\sum_id_i^{-2}b_ib_i^*\big\|^{1/2}\big).$$
Since
 $\wt f\in F^\U\subset (H^c\op H^r)/S^\U$ and $\|\wt f\|=1$,
there exist $\wt f_0\in K_0^c, \wt f_1\in K_1^r$ and $\wt f_2\in
(K^c\op K^r)/\Ga(D)^\perp$ such that
 $$\wt f=\wt f_0+\wt f_1+\wt f_2\quad\mbox{and}\quad
 \|\wt f_j\|\le 1,\quad j=0, 1,2.$$
On the other hand, by the previous decomposition of $(H^c\op
H^r)/S^\U$ we find
 $$\frac{(\el_2(H))^c\op (\el_2(H))^r}{\el_2(S^\U)}=
 (\el_2(K_0))^c\op (\el_2(K_1))^r\op\,
 \frac{(\el_2(K))^c\op (\el_2(K))^r}{\el_2(\Ga(D)^\perp)}\,.$$
Then for any finite sequence $(x_k)\subset S_\8$, by Lemma
\ref{ultrapower} we deduce
 \be
 \big\|\sum_{k} x_k\ot f_k\big\|_{S_\8[F]}
 &\sim_c&
 \max\Big(
 \|\wt f_0\|\,\big\|\sum_kx_k^*x_k\big\|^{1/2},\;
 \|\wt f_1\|\,\big\|\sum_kx_kx_k^*\big\|^{1/2},\\
 && ~~~~~~~~ \big\|\sum_{k} x_k\ot \wt f_2\ot e_k
 \big\|_{S_\8[((\el_2(K))^c\op(\el_2(K))^r)/\el_2(\Ga(D)^\perp)]}\Big).
 \ee
Now write $\wt f_2=(\a_i,\,\b_i)_{i\in I}+\Ga(D)^\perp$ with
 $$\sum_i(|\a_i|^2+|\b_i|^2)\le 2\|\wt f_2\|^2\le2.$$
Thus at most countably many $(\a_i,\,\b_i)$'s are nonzero. On the
other hand, if $\a_i=-d_i\b_i$ for some $i$, then $(\a_i g_i, \b_i
g_i)\in \Ga(D)^\perp$; so this term does not contribute to $\wt
f_2$. Hence, without loss of generality we can assume that the
index set $I$ is equal to $\nz$  and $\a_i\neq -d_i\b_i$ for every
$i$. Let $\wh a=(\wh a_{ik})\in S_\8[(\el_2(K))^c]$ and $\wh
b=(\wh b_{ik})\in S_\8[(\el_2(K))^r]$ such that
 $$\sum_{k} x_k\ot \wt f_2\ot e_k=\id\ot q\ot\id(\wh a, \wh b).$$
Note that $\wh a_{ik}\,,\, \wh b_{ik}\in S_\8$,  the indices $i$
and $k$  correspond to the bases $(g_i)$ of $K$ and $(e_k)$ of
$\el_2$, respectively. Then
 $$\a_ix_k-\wh a_{ik}=-d_i(\b_ix_k-\wh b_{ik}),\quad\forall\; i, k;$$
whence
 $$x_k=\frac{1}{\a_i+d_i\b_i}\, \wh a_{ik}+
 \frac{d_i}{\a_i+d_i\b_i}\,\wh b_{ik}\,{\mathop =^{\rm def}}\,
 a_{ik}+b_{ik}\,.$$
Let
 $$u_c(i)=|\a_i+d_i\b_i|^2\quad\mbox{and}\quad
 u_r(i)=d_i^{-2}|\a_i+d_i\b_i|^2\,.$$
Then
 \be
 \max\Big(\big\|\sum_{i, k}u_c(i)a_{ik}^*a_{ik}
 \big\|^{1/2}\,,\;
 \big\|\sum_{i, k}u_r(i)b_{ik}b_{ik}^*
 \big\|^{1/2}\Big)=
 \max\big(\|\wh a\|_{S_\8[(\el_2(K))^c]}\,,\;
 \|\wh b\|_{S_\8[(\el_2(K))^r]}\big).
 \ee
Therefore, we deduce
 \be
  && \big\|\sum_{k} x_k\ot \wt f_2\ot e_k
 \big\|_{S_\8[((\el_2(K))^c\op(\el_2(K))^r)/\el_2(\Ga(D)^\perp)]}\\
 &&~~ =\inf_{x_k=a_{ik}+b_{ik}}\,
 \max\big(\big\|\big(\sum_{i, k}u_c(i)a_{ik}^*a_{ik}\big)^{1/2}\big\|\,,\;
 \big\|\big(\sum_{i, k}u_r(i)b_{ik}b_{ik}^*\big)^{1/2}\big\|\big).
 \ee
It remains to check the weight condition \eqref{w-con dis}. This
is easy. Indeed,
 \be
 \min(u_c(i),\; u_r(i))\le 2\min(\a_i^2+d_i^2\b_i^2, \;
 d_i^{-2}\a_i^2+\b_i^2)\le 2(\a_i^2+\b_i^2).
 \ee
It follows that
 $$\sum_i \min(u_c(i),\; u_r(i))\le 2\sum_i (\a_i^2+\b_i^2)
 \le 4.$$
Therefore, the theorem is proved. \cqd

\begin{cor}\label{hom 3 spaces}
 Every infinite dimensional space $F\in HQS(C\op R)$
is completely isomorphic to $C$, $R$, $C\cap R$ or $K_{u_c,\,u_r}$
for some positive sequences $u_c$ and $u_r$. Moreover,  in each
case we have $C\cap R\subset F\subset C+R$, up to complete
isomorphism.
 \end{cor}

\pf Let $F$ be represented as in Theorem \ref{hom rep}. Then
 $$\l C\cap\g R\cap K_{u_c,\,u_r}=
 \left\{\begin{array}{ll}
 \displaystyle C\cap R & \textrm{ if } \l>0,\; \g>0\\
 \displaystyle C & \textrm{ if } \l>0,\; \g=0 \\
 \displaystyle R &\textrm{ if } \l=0,\; \g>0\\
 \displaystyle K_{u_c,\,u_r} &\textrm{ if } \l=0,\; \g=0
 \end{array}\right.$$
To show the second part it suffices to prove $C\cap R\subset
K_{u_c,\,u_r}\subset C+R$. Let $(x_k)\subset S_\8$ be a finite
sequence. Consider the decomposition $x_k=a_{ik}+b_{ik}$ given by
 $$a_{ik}=\frac{u_r(i)}{u_c(i)+u_r(i)}\, x_k\quad\mbox{and}\quad
 b_{ik}=\frac{u_c(i)}{u_c(i)+u_r(i)}\, x_k\,.$$
Then
 $$\sum_{i, k}u_c(i)a_{ik}^*a_{ik}
 =\sum_i\frac{u_c(i)\,u_r(i)}{u_c(i)+u_r(i)}\,\sum_kx_k^*x_k\,.
 $$
The first sum on the right hand side is finite by virtue of
\eqref{w-con dis}. We have a similar formula for the row case. It
thus follows that $C\cap R\subset K_{u_c,\,u_r}$. To show the
second inclusion take a decomposition $x_k=a_{ik}+b_{ik}$. Let
 $$\wh a_k=\frac1\mu\,
 \sum_i\frac{u_c(i)\,u_r(i)}{u_c(i)+u_r(i)}\,a_{ik}
 \quad\mbox{and}\quad
 \wh b_k=\frac1\mu\,
 \sum_i\frac{u_c(i)\,u_r(i)}{u_c(i)+u_r(i)}\,b_{ik}\,,$$
where
 $$\mu=\sum_i\frac{u_c(i)\,u_r(i)}{u_c(i)+u_r(i)}\,.$$
Then $x_k=\wh a_k +\wh b_k$. On the other hand, by the convexity
of the operator square function we find
 \be
 \wh a_k^*\wh a_k
 \le\frac1{\mu^2}\,
 \sum_i\frac{u_c(i)\,u_r(i)^2}{\big(u_c(i)+u_r(i)\big)^2}\,
 \sum_iu_c(i)\,a_{ik}^*a_{ik}
 \le \frac1\mu\,\sum_iu_c(i)a_{ik}^*a_{ik}\,.
 \ee
Using similar inequalities for the row terms, we deduce the
inclusion $K_{u_c,\,u_r}\subset C+R$. \cqd

\begin{rk}\label{C+R included}
 Let $u_c$ and $u_r$ be two weights on $(\O,\nu)$
verifying \eqref{w-con}. Using the arguments in the proof of
Corollary~\ref{hom 3 spaces}, one easily checks the following
properties:
 \begin{enumerate}[(i)]
 \item  If
 $$\int_\O u_c\,d\nu<\8\quad\mbox{and}\quad
 \int_\O u_r\,d\nu<\8,$$
then $K_{u_c,\,u_r}=C+R$.
 \item  If $A\subset\O$ is such that
 $$\int_Au_cd\,\nu<\8\quad\mbox{and}\quad
 \int_Au_rd\,\nu<\8,$$
then the parts of $u_c$ and $u_r$ on $A$ do not contribute to
$K_{u_c,\,u_r}$. More precisely, $K_{u_c,\,u_r}$ is completely
isomorphic to $K_{\wt u_c,\, \wt u_r}$, constructed on $A^c$
relative to the weights $\wt u_c=u_c\big|_{A^c}$ and $\wt
u_r=u_r\big|_{A^c}$.
 \end{enumerate}

 \end{rk}


Theorem~\ref{hom rep} can be extended to a more general setting.
Let $X$ and $Y$ be two homogeneous Hilbertian operator spaces. Let
$(u_c,\,u_r)$ be a pair of weights on $(\O,\,\nu)$ satisfying
\eqref{w-con}. We consider the space
 $$G_{u_c,\,u_r}(X,\, Y)=L_2(u_c;X)^X+ L_2(u_r;Y)^Y\,,$$
where $L_2(u_c;X)^X$ denotes $L_2(u_c;X)$ equipped with the
operator space structure of $X$ and similarly for $ L_2(u_r;Y)^Y$.
Let $K_{u_c,\,u_r}(X,\, Y)$ be the constant function subspace of
$G_{u_c,\,u_r}(X,\, Y)$. Reexamining its proof, we find that
Theorem~\ref{hom rep} is still valid with $(C,\, R)$ replaced by
$(X,\,Y)$. Let us record this in the following statement, which is
particularly interesting when $X=C_p$ and $Y=R_p$.

\begin{thm}\label{hom rep bis}
 Let $X$ and $Y$ be two separable infinite dimensional
homogeneous Hilbertian operator spaces. Let $F$ be an infinite
dimensional homogeneous space in $QS(X\op Y)$. Then there exist
two constants $\l, \g\in[0,\;1]$ and two positive sequences
$u_c=(u_c(j))_{j\ge1},\; u_r=(u_r(j))_{j\ge1}$ such that $u_c$ and
$u_r$ satisfy \eqref{w-con dis} and such that $F$ is completely
isomorphic to $\l X\cap\, \g Y\,\cap K_{u_c,\,u_r}(X,\,Y)$.
Moreover, the equivalence constants depend only on the homogeneity
and Hilbertian constants of $X, Y$ and $F$.
 \end{thm}


\section{Homogeneous spaces satisfying a regularity condition}
 \label{Homogeneous spaces satisfying a regularity condition}


The representation of a space $F\in HQS(C\op R)$ given by
Theorem~\ref{hom rep} is far from being unique. The objective of
this section is to give a certain regularity condition on
$(u_c,\,u_r)$ which  ensures that if $F=K_{u_c,\,u_r}$, then up to
equivalence $u_c$ and $u_r$ are uniquely  determined by $F$.  The
first step towards this is the following simple result. In the
sequel $\rz$ will be always equipped with Lebesgue measure.

\begin{prop}\label{reg1}
 Let $(u_c,\,u_r)$ be a pair of weights on $(\O,\,\nu)$
satisfying \eqref{w-con}. Then $K_{u_c,\,u_r}$ is completely
isomorphic to $K_{w_c,\,w_r}$, where $w_c$ and $w_r$ are two
weights on $\rz$ satisfying the following conditions
 \begin{eqnarray}\label{reg-con}
 \left\{\begin{array}{ll}
  w_c=1 \mbox{ on } \rz_+\,,
 w_c \mbox{ is nondecreasing and left continuous};\\
 w_r=1 \mbox{ on } \rz_-\,,
 w_c \mbox{ is nonincreasing and right continuous};\\
  \displaystyle\int_{-\8}^0w_c(s)\,ds=1
 =\int_0^\8w_r(s)\,ds.
 \end{array} \right.
 \end{eqnarray}
 \end{prop}

\pf Consider the following partition of $\O$:
 $$A_k=\big\{\o\in\O\,:\, 2^{-k-1}\le\frac{u_r(\o)}{u_c(\o)}
 <2^{-k}\big\},\quad k\in\zz.$$
Let
 $$\wt u_c(k)=u_c(A_k)
 \quad\mbox{and}\quad
 \wt u_r(k)=2^{-k}u_c(A_k)\,\quad k\in\zz,$$
where
 $u_c(A)=\int_Au_c\,d\,\nu$
for $A\subset\O$. By Proposition~\ref{cont-dis}, $K_{u_c,\,u_r}$
is completely isomorphic to $K_{\wt u_c,\,\wt u_r}$. Now put
 $$s_0=0,\quad s_{-1}=-1\quad\mbox{and}\quad
 s_k=-u_r(A_{k+1})-u_r(A_{k+2})-\cdots-u_r(A_{-1})-1
 \mbox{ for } k\le -2.$$
 Define
 $$w_c(s)=\frac1{2\l}\,\sum_{k\le-2}2^{k}\un_{(s_k,\;s_{k+1}]}(s)\;
 \mbox{ for }\; s\in(-\8,\;-1].$$
where
 $$\l=\max\big(1,\;\int_\O\min(u_c,\;u_r)d\,\nu\big).$$
Then
 \be
 \int_{-\8}^{-1}w_c(s)\,ds
 &=&\frac1{2\l}\,\sum_{k\le-2}2^{k}u_r(A_{k+1})\\
 &\le& \frac1{2\l}\,\sum_{k\le-2}\int_{A_{k+1}}u_cd\,\nu
 \le\frac1{2\l}\,\int_\O\min(u_c,\;u_r)d\,\nu
 \le \frac12\,.
 \ee
Completing $w_c$ by
 $$w_c(s)=1-\int_{-\8}^{-1}w_c(s)\,ds\;
 \mbox{ for } s\in(-1,\;0] \quad\mbox{and}\quad
 w_c(s)=1 \mbox{ for } s>0,$$
we see that $w_c$ is nondecreasing and
 $$\int_{-\8}^0w_c(s)\,ds=1.$$
The second function $w_r$ is defined similarly by using $u_r$
instead of $u_c$. Indeed, letting
$$s_{1}=1\quad\mbox{and}\quad
 s_k=1+u_c(A_{1})+\cdots +u_c(A_{k-1})
 \mbox{ for } k\ge 2,$$
we define $w_r$ by
 $$w_r=\un_{(-\8,\; 0)} +
 \big(1-\int^{\8}_{1}w_r(s)ds\big)\,\un_{[0,\; 1)}+
 \frac1{2\l}\,\sum_{k\ge1}2^{-k-1}\un_{[s_k,\;s_{k+1})}.$$
Then by Proposition~\ref{cont-dis} (and Remark~\ref{C+R
included}), $K_{w_c,\,w_r}$ is completely isomorphic to $K_{\wt
u_c,\,\wt u_c}$, so to $K_{u_c,\,u_r}$ too. \cqd

\begin{rk}\label{reg1 equiv}
 The last condition in \eqref{reg-con} is not essential
for $K_{w_c,\,w_r}$. What really matters is the behavior of the
weights at infinity. More precisely, if $v_c$ and $v_r$ are two
other locally integrable weights on $\rz$ such that
 $$c_1 v_c(s)\le w_c(s)\le c_2 v_c(s)\quad\mbox{and}\quad
 c_1 v_r(s)\le w_r(s)\le c_2 v_r(s)$$
for two positive constants $c_1$, $c_2$ and for all $s\in\rz$ with
$|s|$ sufficiently big, then $K_{w_c,\,w_r}$ and $K_{v_c,\,v_r}$
are completely isomorphic. In particular, if $w_c$ and $w_r$
satisfy the two first conditions in \eqref{reg-con} and are
integrable on $\rz_-$ and $\rz_+$, respectively, then we can
adjust them, up to equivalence, so that they satisfy the last
condition too. This is achieved by the arguments in the proof of
Proposition~\ref{reg1} and Remark~\ref{C+R included}.
 \end{rk}

In the following we will often assume that $w_c$ and $w_r$ satisfy
\eqref{reg-con}. Note that the left (right) continuity is assumed
only for convenience. In fact, by perturbation and
Remark~\ref{reg1 equiv} we can even assume that $w_c$ and $w_r$
are continuous. For the same reason, we can also assume that $w_c$
(resp. $w_r$) is increasing on $\rz_-$ (resp. decreasing on
$\rz_+$). On the other hand, if $w_r=0$ on $[s_0,\; \8)$ for some
$s_0\ge0$, then by Remarks~\ref{positive weight} and \ref{C+R
included} we need only to consider the part on $\rz_-$. Thus for
presentation simplicity we will assume that both $w_c$ and $w_r$
are strictly positive on $\rz$ in the sequel.

We will show that under a certain regularity condition $w_c$ and
$w_r$ are uniquely determined, up to equivalence at infinity,  by
the operator space structure of $K_{w_c,\, w_r}$, i.e., by the
fundamental functions of $K_{w_c,\, w_r}$. To this end we require
some preparations and two auxiliary pairs of functions associated
with $(w_c,\, w_r)$. We will concentrate only on $w_r$ since $w_c$
can be dealt with by using the symmetry $t\mapsto -t$.

\medskip

Let $w$ be an integrable nonincreasing right continuous positive
function on $\rz_+$. Consider the functions $h$ and $g$ on $\rz_+$
defined by
 \beq\label{def hg}
 h(t)=\int_t^\8w(s)ds\quad\mbox{and}\quad
 g(s)=\inf_{t>0}\max\big(t,\; sh(t)\big).
 \eeq
The following is entirely elementary.

\begin{lem}\label{reg2}
 Let $w,\,h,\,g$ be as above.
 \begin{enumerate}[\rm(i)]
 \item $h$ is decreasing, continuous and derivable
from right with $-w$ as right derivative.
 \item $g$ is increasing, continuous and $g(s)\to\8$
as $s\to\8$.
 \item For any $s\in\rz_+$ the infimum defining $g(s)$ is attained
at a unique $t$ such that $t=sh(t)$. Moreover $($with $g^{-1}$
denoting the inverse of $g)$,
 \beq\label{g-h}
 \frac{h(t)}t=\frac1{g^{-1}(t)}\,,\quad\forall\; t>0.
 \eeq
 \end{enumerate}
 \end{lem}

\pf (i) and (ii) are obvious except the injectivity of $g$ which
will follow from (iii). To prove (iii) we observe that by the
continuity of $h$ and the fact that $h(t)\to0$ as $t\to\8$, for
any $s\in\rz_+$ there exists $t$ such that $g(s)=\max\big(t,\;
sh(t)\big)$. If $t<sh(t)$, then by the continuity of $h$ once more
there exists $t_1$ such that $t_1>t$ and $t_1<sh(t_1)$; whence
 $$g(s)\le\max\big(t_1,\; sh(t_1)\big)=sh(t_1)<sh(t)=
 \max\big(t,\; sh(t)\big)=g(s),$$
which is a contradiction. Similarly, $t>sh(t)$ cannot happen.
Therefore, $t=sh(t)$; so $g(s)=t$ and $g(s)=sh(t)$. This implies
in particular that $g$ is injective, so increasing. Consequently,
its inverse $g^{-1}$ exists and $s=g^{-1}(t)$. Together with
$t=sh(t)$, this yields \eqref{g-h}. \cqd

\medskip

We will employ Orlicz space techniques, one of the main novelties
of the paper. Recall that an Orlicz function is a convex function
$\f$ on $[0,\;\8)$ with $\f(0)=0$. The Orlicz sequence space
$\el_\f$ consists of all complex sequences $x=(x_n)$ such that
 $$\sum_n\f\big(\frac{|x_n|}{\l}\big)<\8$$
for some $\l>0$  and is equipped with the norm
 $$\|x\|_\f=\inf\big\{\l>0\;:\;
 \sum_n\f\big(\frac{|x_n|}{\l}\big)\le1\big\}.$$
Note that $\el_\f$ depends, up to an equivalent norm, only on the
values of $\f$ in a neighborhood of $0$. It is well-known that an
Orlicz function $\f$ is uniquely determined, up to equivalence at
$0$, by its fundamental sequence $(\f_n)_{n\ge1}$, where
 $$\f_n=\big\|\sum_{k\le n}e_k\big\|_\f=\frac1{\f^{-1}(\frac1n)}\,.$$
This is equivalent to saying that up to an equivalent norm,
$\el_\f$ is uniquely determined by $(\f_n)$ (see
\cite[Proposition~4.a.5]{LT-I}). Here $\f^{-1}$ denotes the
generalized inverse of $\f$, i.e.,
$\f^{-1}(t)=\sup\big\{s\,:\,\f(s)\le t\big\}$. Note that $\f$ is
invertible except in the trivial case where $\f$ vanishes in a
neighborhood of $0$.

\medskip

Now let $(w_c,\,w_r)$ be a pair of positive weights on $\rz$
satisfying \eqref{reg-con}.  The functions $h$ and $g$ associated
to $w_r$ as in \eqref{def hg} will be denoted by $h_r$ and $g_r$,
respectively. Accordingly, the two functions on $\rz_+$ associated
to the function $t\mapsto w_c(-t)$ will be denoted by $h_c$ and
$g_c$. Define $\f_c,\,\f_r:\,\rz_+\to\rz_+$ by
 $$\f_c(t)=t^2\,h_c(\frac1{t^2})\quad\mbox{and}\quad
 \f_r(t)=t^2\,h_r(\frac1{t^2})$$
(with $\f_c(0)=\f_r(0)=0$). It is easy to check that $\f_c$ and
$\f_r$ are convex functions. Indeed, using the equality
 $$\f_r(t)=t^2\int_{t^{-2}}^\8w_r(s)ds, $$
we find
 $$\f_r''(t)=2\int_{t^{-2}}^\8w_r(s)ds +\frac2{t^2}\,
 w_r(\frac1{t^2}) -\frac4{t^4}\,w'_r(\frac1{t^2}),$$
where the derivatives are right derivatives. Since $w_r$ is
nonincreasing, $\f_r''(t)\ge0$, so $\f_r$ is convex.

\medskip
The next lemma contains the main idea of our uniqueness theorem on
$K_{w_c,\,w_r}$ and also explains the reason for the introduction
of $h_r$ and $g_r$. Its proof uses an Orlicz space argument which
will also play a key role in the proof of Theorem~\ref{c1s orlicz}
below. We will consider only the diagonal elements of
$R[K_{w_c,\,w_r}]$ although we can determine the norm of any
element (see the proof of Theorem~\ref{c1s orlicz} for more
details). Recall that $(e_k)$ denotes the canonical basis of
$K_{w_c,\,w_r}$ and $(e_{k1})$ (resp. $(e_{1k})$) that of $C$
(resp. $R$).

\begin{lem}\label{orlicz cr}
 For any finite sequence $(x_k)\subset\cz$ we
have
 \beq\label{o-norm}
 \big\|\sum_kx_ke_{1k}\ot e_k\big\|_{R[K_{w_c,\,w_r}]}
 \sim_c\|(x_k)\|_{\f_r}.
 \eeq
Moreover, for any $n\in\nz$
 \beq\label{o-norm1}
 \big\|\sum_{k=1}^ne_{1k}\ot e_k\big\|^2_{R[K_{w_c,\,w_r}]}
 \sim_c  g_r(n).
 \eeq
A similar statement holds for the column case.
 \end{lem}

\pf We consider only the row case, the column part being treated
similarly. We will use the following auxiliary function
 $$\mathsf h_r(t)=\inf\big\{w_r(A^c)\,:\, A\subset\rz,\;
 w_c(A)\le t\big\},\quad t>0.$$
We claim that
 $$h_r(t)\le \mathsf h_r(t)\quad\mbox{and}\quad
 \mathsf h_r(1+t)\le h_r(t).$$
Indeed, first observe that
 $$\sup\big\{w_r(A)\,:\, A\subset\rz_+\,,\; |A|\le t\big\}=\int_0^t
 w_r(s)ds$$
for $w_r$ is nonincreasing, where $|A|$ denotes the Lebesgue
measure of $A$. Then we find
 \be
 h_r(t)
 &=&\int_t^\8 w_r(s)ds=1-\int_0^t w_r(s)ds\\
 &=&\inf\big\{w_r(A^c\cap\rz_+)\,:\, A\subset\rz\,,\;
 |A\cap\rz_+|\le t\big\}\\
 &=&\inf\big\{w_r(A^c\cap\rz_+)\,:\, A\subset\rz\,,\;
 w_c(A\cap\rz_+)\le t\big\}.
 \ee
Therefore, $h_r\le \mathsf h_r$. Since $w_c(\rz_-)=1$, we also
have $ \mathsf h_r(1+t)\le h_r(t)$.

Let $\phi_r(t)=t^2\,\mathsf h_r(t^{-2})$. Then the preceding claim
shows that $\f_r$ and $\phi_r$ are equivalent at $0$ with
universal equivalence constants. Thus we need only to prove
\eqref{o-norm} with $\phi_r$ in place of $\f_r$. To this end first
recall the following elementary identifications
 $$R[\el_2^c]=R[C]=\mathbb K(\el_2)\quad\mbox{and}\quad
 R[\el_2^r]=R[R]=(S_2)^r\,,$$
where $\mathbb K(\el_2)$ denotes the space of compact operators on
$\el_2$. Then we find
 $$R[G_{w_c,\,w_r}]=R[L_2^c(w_c;\el_2)]+R[L_2^r(w_r;\el_2)]
 =L_2^c(w_c)\ot_{\min}\mathbb K(\el_2) + L_2^r(w_r;S_2)\,.$$
Since we are dealing only with diagonal elements of
$R[K_{w_c,\,w_r}]$, we use the diagonal projection $D$ on $\mathbb
K(\el_2)$ and $S_2$. The complete contractivity of $D$ on $\mathbb
K(\el_2)$ implies that $\id\ot D$ is a contraction on
$L_2^c(w_c)\ot_{\min}\mathbb K(\el_2)$. On the other hand, it is
trivial that $\id\ot D$ is also a contraction on $L_2^r(w_r;S_2)$.
Therefore, $\id\ot D$ is a contractive projection from
$R[G_{w_c,\,w_r}]$ onto its diagonal part, which is
$L_2^c(w_c)\ot_{\min}c_0 + \big(\el_2(L_2(w_r)\big)^r$. Since
$R[K_{w_c,\,w_r}]$ is the constant function subspace of
$R[G_{w_c,\,w_r}]$, for any finite sequence $x=(x_k)\subset\cz$ we
then deduce
 \be
 \big\|\sum_kx_ke_{1k}\ot e_k\big\|_{R[K_{w_c,\,w_r}]}
 =\inf\big\{\sup_k\big\|a_k\big\|_{L_2(w_c)}
 +\big(\sum_k\big\|b_k\big\|_{L_2(w_r)}^2\big)^{1/2}\big\},
 \ee
where the infimum runs over all decompositions $x_k=a_k+b_k$ a.e.
on $\rz$ with $a_k\in L_2(w_c)$ and $b_k\in L_2(w_r)$.  Put
 $$\tnorm{x}=\inf\big\{
 \sup_k|x_k|\,w_c(A_k)^{1/2}+
 \big(\sum_k|x_k|^2w_r(A_k^c)\big)^{1/2}\,:\,A_k\subset\rz\big\}.$$
It is then easy to check that
 $$\frac12\,\tnorm{x}\le
 \big\|\sum_kx_ke_{1k}\ot e_k\big\|_{R[K_{w_c,\,w_r}]}
 \le \tnorm{x}.$$
Indeed, given a decomposition $x_k=a_k+b_k$ as above, define
 $$A_k=\big\{s\in\rz\,:\, |a_k(s)|\ge\frac{|x_k|}2\big\}.$$
Then $|b_k|>|x_k|/2$ on $A_k^c$. Thus
 $$\frac12\,\tnorm{x}\le
 \big\|\sum_kx_ke_{1k}\ot e_k\big\|_{R[K_{w_c,\,w_r}]}\,.$$
The inverse inequality follows by taking the special decomposition
$x_k=x_k\un_{A_k} + x_k\un_{A_k^c}$.

Now assume $\big\|\sum_{k=1}^nx_ke_{1k}\ot
e_k\big\|_{R[K_{w_c,\,w_r}]}<1/2$. Without loss of generality,
assume that $x_k\neq0$ for every $1\le k\le n$. Then we find
$A_k\subset\rz$ such that
 $$\sup_k|x_k|\,w_c(A_k)^{1/2}
 +\big(\sum_k|x_k|^2w_r(A_k^c)\big)^{1/2}<1.$$
Therefore,
 $$w_c(A_k)\le \frac1{|x_k|^{-2}}\,,\quad\forall\; 1\le k\le n$$
and
 $$\sum_{k=1}^n\phi_r(|x_k|)=
 \sum_{k=1}^n|x_k|^2\,\mathsf h_r(\frac1{|x_k|^{-2}})
 \le \sum_k|x_k|^2\, w_r(A_k^c)\le 1.$$
Thus $\|(x_k)\|_{\phi_r}\le1$. Conversely, assume that
$\|(x_k)\|_{\phi_r}<1$. Then there exist $A_k\subset\rz$ such that
 $$w_c(A_k)\le \frac1{|x_k|^{-2}}\quad\mbox{and}\quad
 \sum_k|x_k|^2 w_r(A_k^c)\le1.$$
It follows that $\tnorm{x}\le 2$; so
$\big\|\sum_{k=1}^nx_ke_{1k}\ot e_k\big\|_{R[K_{w_c,\,w_r}]}\le
2$. Therefore, \eqref{o-norm} is proved. On the other hand,
\eqref{g-h} implies
 $$\frac1{\sqrt{g_r(s)}}=\f^{-1}(\frac1s)\,,\quad\forall\; s>0.$$
\eqref{o-norm1} is then a particular case of \eqref{o-norm} since
 $$
  \big\|\sum_{k=1}^ne_k\big\|_{\f_r}=\frac1{\f_r^{-1}(\frac1n)}
  =\sqrt{g_r(n)}\,.
  $$
\cqd

\medskip

Let $F=K_{w_c,\,w_r}$. Then  \eqref{o-norm1} shows that $R[F]$ is
determined (up to equivalence) by $g_r$, so by $w_r$ too. In order
to reverse this procedure, i.e., to determine $w_r$ by $R[F]$, we
must know how to recover $w_r$ from $g_r$. This is possible under
a certain regularity condition on $w_r$ (or $g_r$), as shown by
the following lemma.

\begin{lem}\label{reg3}
 Let $w,\,h,\,g$ be as in Lemma~\ref{reg2}.
 Then the following properties are equivalent:
 \begin{enumerate}[\rm(i)]
 \item There exist positive constants $c_1,\, d_1$ and
$\a_1$ with $\a_1>1$ such that
 \beq\label{reg w}
 c_1\left(\frac{t}{s}\right)^{\a_1}\le\frac{w(s)}{w(t)}
 \quad\mbox{and}\quad
 \frac{w(s)}{w(2s)}\le d_1\,,\quad\forall\;t\ge s\ge1.
 \eeq
 \item There exist positive constants $c_2,\, d_2$
and $\a_2$ such that
 \beq\label{reg h}
 c_2\left(\frac{t}{s}\right)^{\a_2}\le\frac{h(s)}{h(t)}
 \quad\mbox{and}\quad
 \frac{h(s)}{h(2s)}\le d_2\,,\quad\forall\;t\ge s\ge1.
 \eeq
\item There exist positive constants $c,\, d$ and $\b$ with
$0<\b<1$ such that
 \beq\label{reg g}
 \frac{g(t)}{g(s)}\le d\left(\frac{t}{s}\right)^{\b}
 \quad\mbox{and}\quad
 \frac{g^{-1}(2s)}{g^{-1}(s)}\le c\,,\quad\forall\;t\ge s\ge1.
 \eeq
 \end{enumerate}
Moreover, if one of  {\rm (i)-(iii)} is satisfied, then
 \beq\label{g-w}
 c'\,w(t)\le \frac1{g^{-1}(t)}\le d'\,w(t),\quad\forall\;t\ge1\,,
 \eeq
where the constants $c'$ and $d'$ depend only on the relative
constants in {\rm (i)-(iii)}.
 \end{lem}

\pf First note that by \eqref{g-h} for $s>0$ and $t=g(s)$ we have
 $$h(t)=\frac{g(s)}s=\frac{t}{g^{-1}(t)}\,.$$
Using this we easily check  (ii) $\Leftrightarrow$ (iii) (with
$\b=(1+\a_2)^{-1}$) and thus omit the details. Now assume (ii).
Let us show \eqref{g-w} and (i). We have
 $$w(t)\le 2\,\frac{h(t/2)}t\le 2d_2\frac{h(t)}t\,,
 \quad \forall\;t\ge2\,.$$
On the other hand, since $w$ is decreasing, $h$ is convex. Thus
for any $0<t<s<r$ we have
 $$\frac{h(s)-h(t)}{s-t}\le \frac{h(r)-h(t)}{r-t}\,.$$
Fixing $r$ and letting $s\to t$, we deduce
 $$-w(t)\le \frac{h(r)-h(t)}{r-t}\,.$$
Now choose $r=\l t$ with $c_2\l^{\a_2}\ge 2$. Then
 $$h(\l t)\le \frac{h(t)}{c_2\l^{\a_2}}\le \frac{h(t)}2\,.$$
It then follows that
 $$w(t)\ge \frac1{2(\l -1)}\, \frac{h(t)}t\,.$$
Therefore,
 $$w(t)\sim \frac{h(t)}t=\frac1{g^{-1}(t)}\,.$$
This is \eqref{g-w} and also yields (i) with $\a_1=1+\a_2$. It
remains to show (i) $\Rightarrow$ (ii). Assuming (i), we have
 $$h(t)\ge tw(2t)\ge d_1^{-1}\,tw(t),\quad\forall\;t\ge1\,.$$
For the converse inequality we find
 $$h(t)\le c_1^{-1}t^{\a_1}w(t)\int_t^\8s^{-\a_1}ds
 =c_1^{-1}(\a_1-1)^{-1}t w(t).$$
Combining this with the previous inequality we get $h(t)\sim
tw(t)$. We then deduce (ii) from (i). \cqd

\begin{rk}
 The second inequality in \eqref{reg w} is satisfied if
there exist  positive constants $d_1$ and $\b_1$ with $\b_1>1$
such that
 \beq\label{reg wbis}
 \frac{w(s)}{w(t)}\le d_1\left(\frac{t}{s}\right)^{\b_1}\,,
 \quad\forall\; t\ge s\ge1.
 \eeq
A similar remark applies to \eqref{reg h} and \eqref{reg g}.
However, in the case of \eqref{reg g} we can formulate a condition
directly on $g$ as follows: there exist positive constants $c$ and
$\a$ with $0<\a<1$ such that
 \beq\label{reg gbis}
 c\left(\frac{y}{x}\right)^\a\le\frac{g(y)}{g(x)}\,,
 \quad\forall\; y\ge x\ge1.
 \eeq
 \end{rk}

\begin{rk}\label{reg2 inverse}
 Lemma~\ref{reg2} defines the function $g$ starting from
$w$. Under the regularity conditions in Lemma~\ref{reg3} we are
able to reverse this procedure. Namely, let $\mathsf g$ be an
increasing continuous function on $\rz_+$ satisfying \eqref{reg
g}. Let $w=1/\mathsf g^{-1}$. Let $g$ be the function associated
to $w$ by \eqref{def hg}. Then by Lemma~\ref{reg3} $g$ is
equivalent to $\mathsf g$ at infinity: there exist two positive
constants $c$ and $d$ such that
 $$cg(x)\le \mathsf g(x)\le d\, g(x),\quad\forall\; x\ge1.$$
 \end{rk}

We are finally ready to prove our uniqueness theorem on
$K_{w_c,\,w_r}$.

\begin{thm}\label{uni}
 Let $\mathsf g_c$ and $\mathsf g_r$ be two
increasing continuous positive functions on $\rz_+$ both
satisfying the regularity condition \eqref{reg g} for some
positive constants $c,\, d$ and $\b$ with $0<\b<1$. Then there
exists, up to complete isomorphism, exactly one  space $F\in
HQS(C\op R)$ such that
 \beq\label{reg fund1}
 \big\|\sum_{k=1}^ne_{k1}\ot e_k\big\|_{C[F]}^2
 \sim \mathsf g_c(n),\quad
 \big\|\sum_{k=1}^ne_{1k}\ot e_k\big\|_{R[F]}^2
 \sim \mathsf g_r(n)\quad\mbox{uniformly in }\; n\in\nz,
 \eeq
where $(e_k)$ is an orthonormal basis of $F$. Moreover, $F$ can be
chosen to be $K_{w_c,\, w_r}$ with a pair $(w_c,\,w_r)$ of weights
on $\rz$ satisfying \eqref{reg-con} and such that
 $$w_c(-t)\sim \frac1{\mathsf g_c^{-1}(t)}
 \quad\mbox{and}\quad
 w_r(t)\sim \frac1{\mathsf g_r^{-1}(t)}\quad\mbox{for}\quad t\ge1.$$
 \end{thm}

\pf By Remarks~\ref{reg2 inverse} and \ref{reg1 equiv} we may
perturb $(\mathsf g_c^{-1}\,, \mathsf g_r^{-1})$ into a pair
$(w_c,\,w_r)$ of weights on $\rz$ verifying \eqref{reg-con} and
\eqref{reg w}. Therefore, thanks to Lemma~\ref{orlicz cr}, the
space $F=K_{w_c,\, w_r}$ satisfies \eqref{reg fund1}. Now let
$F\in HQS(C\op R)$ satisfy \eqref{reg fund1}. Using
Theorem~\ref{hom rep} and Proposition~\ref{reg1}, we can assume
$F=\l C\cap\g R\cap K_{w_c,\,w_r}$. Then we must have $\l=\g=0$.
Indeed, if, for instance, $\l>0$ and $\g=0$, then $F=C$ and so
$\|\sum_{k=1}^ne_{k1}\ot e_k\big\|_{C[F]}=\sqrt n\,$, a contradiction with
\eqref{reg fund1} and the regularity of $\mathsf g_c$. Therefore,
$F=K_{w_c,\,w_r}$. Let $g_c$ and $g_r$ be the functions associated
to $w_c$ and $w_r$ as in \eqref{def hg}. Then by Lemma~\ref{orlicz
cr}, $g_c\sim \mathsf g_c$ and $g_r\sim \mathsf g_r$. Therefore,
$g_c$ and $g_r$ satisfy \eqref{reg g} too. Thus by \eqref{g-w} we
deduce
 $$w_c(-t)\sim \frac1{g_c^{-1}(t)}\sim \frac1{\mathsf g_c^{-1}(t)}
 \quad\mbox{and}\quad
 w_r(t)\sim \frac1{g_r^{-1}(t)}\sim \frac1{\mathsf g_r^{-1}(t)}
 \quad\mbox{for}\quad t\ge1.$$
Hence $F$ is completely isomorphic to the space constructed
previously. This shows the uniqueness part. \cqd

\medskip

The preceding results suggest the following

\begin{definition}
 Let $F$ be a homogenous Hilbertian operator space. Let
 $$\F_{c, F}(n)
 =\big\|\sum_{k=1}^ne_{k1}\ot e_k\big\|_{C[F]}^2
 \quad\mbox{and}\quad
 \F_{r, F}(n)=
 \big\|\sum_{k=1}^ne_{1k}\ot e_k\big\|_{R[F]}^2\,,\quad n\in\nz,$$
where $(e_k)$ is a basis of $F$ which is the image of an
orthonormal basis under an isomorphism between $F$ and a Hilbert
space. We call $\F_c$ and $\F_r$ the \emph{fundamental functions}
of $F$. $F$ is said to be \emph{regular} if there exist positive
constants $c,\,d$ and $\a,\,\b$ with $0<\a\le \b<1$ such that
 \be
 c\left(\frac{n}{k}\right)^\a\le
 \frac{\F_{c, F}(n)}{\F_{c, F}(k)}
 \le d\left(\frac{n}{k}\right)^\b\quad\mbox{and}\quad
 c\left(\frac{n}{k}\right)^\a\le
 \frac{\F_{r, F}(n)}{\F_{r, F}(k)}
 \le d\left(\frac{n}{k}\right)^\b\,,\quad\forall\; n\ge k\ge1.
 \ee
 \end{definition}

We will denote $\F_{c, F}$ and $\F_{r, F}$ simply by $\F_{c}$ and
$\F_{r}$, respectively whenever no confusion can occur. In the
sequel, these functions will be systematically extended to
piecewise linear functions on $\rz_+$ with $\F_{c}(0)=\F_{c}(1)$
and $\F_{r}(0)=\F_{r}(1)$. In fact, we can even assume, by
perturbation, that they are increasing on $\rz_+$.

\begin{rk}\label{uni bis}
Note that the homogeneity of $F$ implies
 $$\frac{\F_c(n)}{\F_c(k)}
 \le d\,\frac{n}{k}\quad\mbox{and}\quad
 \frac{\F_r(n)}{\F_r(k)}
 \le d\,\frac{n}{k}\,.$$
Thus the regularity of $F$ does not require much more than
necessary.
 \end{rk}

\begin{thm}\label{dual}
 Let $F$ be a $1$-homogeneous $1$-Hilbertian operator space.
Then
 $$\F_{c, F}(n)\F_{c, F^*}(n)=n\quad\mbox{and}\quad
 \F_{r, F}(n)\F_{r, F^*}(n)=n, \quad\forall\;n\in\nz.$$
 \end{thm}

\pf By symmetry, it suffices to deal with the column case.  Let
$n\in\nz$. Considering an $n$-dimensional subspace instead of $F$,
we can assume $\dim F=n$. We will use the duality equality:
$C^n[F]^*=C^n_1[F^*].$  Fix an orthonormal basis $(e_k)$ of $F$.
We consider, as usual, the elements $x$ in $C^n[F]$ or
$C^n_1[F^*]$ as matrices with respect to the bases $(e_{k1})$ and
$(e_k)$:
 $$x=\sum_{i,j=1}^nx_{ij}e_{i1}\ot e_j\,.$$
By a slight abuse of notation, we use $\id$ to denote the identity
matrix as well as the corresponding identity elements in these
spaces. Then
 $$\|\id\|_{C^n[F]}=
 \sup\big\{\frac{|\Tr(x)|}{\|x\|_{C^n_1[F^*]}}\,:\,
 x\in C^n_1[F^*]\big\},$$
where $\Tr$ denotes the usual matrix trace. We claim that the
supremum above is attained at $\id$. Indeed, let $\mathcal U$
denote  the group of $n\times n$ unitary matrices, equipped with
Haar measure. For any $u\in\mathcal U$ and $x\in C^n_1[F^*]$ we
have
 $$\Tr(uxu^*)=\Tr(x)\quad\mbox{and}\quad
 \|uxu^*\|_{C^n_1[F^*]}=\|x\|_{C^n_1[F^*]}\,,$$
where the second equality follows from homogeneity. Let
 $$\wt x=\int_{\mathcal U}uxu^*du.$$
Then
 $$\wt x=\frac{\Tr(x)}n\,\id\quad\mbox{and}\quad
 \|\wt x\|_{C^n_1[F^*]}\le\|x\|_{C^n_1[F^*]}\,.$$
Thus the claim follows. Our next task is to calculate
$\|\id\|_{C^n_1[F^*]}$. To this end recall that by
\cite[Theorem~1.5]{pis-ast} we have
 $$\|\id\|_{C^n_1[F^*]}=\inf\big\{
 \|a\|_{S_2^n}\|y\|_{C^n[F^*]}\,:\, \id=ay,\;
 a\in S_2^n,\; y\in C^n[F^*]\big\}.$$
We will show that the infimum is attained at the trivial
factorization with $a=\id$ and $y=\id$. Let $\id=ay$ be a
factorization as above. Then $a=y^{-1}$. Let $u, v\in\mathcal U$
such that $uyv=\d$ is a diagonal matrix with positive diagonal
entries. By the homogeneity of $F, C^n$ and the unitary invariance
of $S_2^n$ we have
 $$\|\d\|_{C^n[F^*]}=\|y\|_{C^n[F^*]}\quad\mbox{and}\quad
 \|\d^{-1}\|_{S_2^n}=\|v^*au^*\|_{S_2^n}= \|a\|_{S_2^n}\,.$$
Thus we can replace $y$ by $\d$ and $a$ by $\d^{-1}$. Then as
before, we have
 $$\int_{\mathcal U}u \d u^* du=\frac{\d_1+\cdots + \d_n}n\,\id,$$
where $\d_1,..., \d_n$ are the diagonal entries of $\d$. Again by
homogeneity, we find
 $$\frac{\d_1+\cdots +  \d_n}n\,\|\id\|_{C^n[F^*]}
 \le\|\d\|_{C^n[F^*]}\,.$$
Similarly,
 $$\frac{\d_1^{-1}+\cdots  + \d_n^{-1}}n\,\|\id\|_{S_2^n}
 \le\|\d^{-1}\|_{S_2^n}\,.$$
We then deduce, by the Cauchy-Schwarz inequality, that
 \be
 \sqrt n\,\|\id\|_{C^n[F^*]}
 &=&\|\id\|_{S_2^n}\|\id\|_{C^n[F^*]}\\
 &\le& \frac{\d_1^{-1}+\cdots  +
 \d_n^{-1}}n\,\frac{\d_1+\cdots  + \d_n}n\,
 \|\id\|_{S_2^n}\|\id\|_{C^n[F^*]}\\
 &\le& \|\d^{-1}\|_{S_2^n}\|\d\|_{C^n[F^*]}\,.
 \ee
This yields the announced assertion. Therefore,
 $$\|\id\|_{C^n_1[F^*]}=\sqrt n\,\|\id\|_{C^n[F^*]}\,.$$
Combining the previous two parts, we finally get
 $$\|\id\|_{C^n[F]}
 =\frac{\Tr(\id)}{\sqrt n\,\|\id\|_{C^n[F^*]}}
 =\frac{\sqrt n}{\|\id\|_{C^n[F^*]}}\,.$$
Recalling that $\F_{c, F}=\|\id\|_{C^n[F]}^2$ and $\F_{c,
F^*}=\|\id\|_{C^n[F^*]}^2$, we then deduce the theorem.\cqd

\begin{cor}\label{reg-duality}
 Let $F$ be a homogeneous Hilbertian space. Then $F$
is regular iff $F^*$ is regular.
 \end{cor}


\section{Completely 1-summing maps}
 \label{Completely 1-summing maps}


In this section we deal  with completely $1$-summing  maps between
two homogeneous spaces $E$ and $F$ in $QS(C\op R)$. We will show
that $\Pi_1^o(E,\, F)$ coincides with a  Schatten-Orlicz class
$S_\f$. We first recall the relevant definition. Let $\Phi$ be a
symmetric sequence space on $\nz$ in the sense of \cite{LT-I}.
Recall that the unitary ideal $S_\Phi$ is the space of all compact
operators $x$ on $\el_2$ such that $(s_n(x))_{n\ge1}$ belongs to
$\Phi$. Here $(s_n(x))_{n\ge1}$ denotes the sequence of singular
values of $x$, i.e., $(s_n(x))_{n\ge1}$ is the sequence of the
eigenvalues of $|x|=\sqrt{x^*x}$ ranged in decreasing order and
repeated according to multiplicity. $S_\Phi$ is equipped with the
norm
 $$\|x\|_{\Phi}=\|(s_n(x))_{n\ge1}\|_\Phi\,.$$
Note that $S_\Phi$ is an ideal of $B(\el_2)$ and its norm is
unitarily invariant. Namely, $\|ux\|_\Phi=\|xu\|_\Phi=\|x\|_\Phi$
for any unitary $u\in B(\el_2)$ and any $x\in S_\Phi$. Conversely,
any unitarily invariant norm on an ideal of compact operators on
$\el_2$ is of this form. We refer to \cite{GK-Sp} and \cite{simon}
for more information.  We are particularly interested in the case
where $\Phi$ is an Orlicz sequence space $\el_\f$. In this case
the unitary ideal $S_{\el_\f}$ is simply denoted by $S_\f$. In
particular, $S_\f=S_p$ if $\f(t)=t^p$ with $p<\8$.  Let us adopt
the following convention. Let $E$ and $F$ be two (infinite
dimensional) homogeneous Hilbertian operator spaces. Identifying
both $E$ and $F$ with $\el_2$ isomorphically, we will view maps
from $E$ to $F$ as operators on $\el_2$.

\begin{rk}\label{c1s ideal}
 Let $E$ and $F$ be two $1$-homogenous $1$-Hilbertian
operator spaces. Then the norm of $\Pi_1^o(E,\, F)$ is unitarily
invariant. This reduces the determination of the whole space
$\Pi_1^o(E,\, F)$ to that of the subspace of diagonal operators
with respect to two given orthonormal bases in $E$ and $F$,
respectively. Note, however, that $\Pi_1^o(E,\, F)$ may contain
non compact operators in general.
 \end{rk}

Recall that an Orlicz function $\f$ satisfies the
$\Delta_2$-condition (at $0$) if there exists a positive constant
$\l$ such that $\f(2t)\le\l\,\f(t)$ for all $t>0$ (in a
neighborhood of $0$). It is well-known that $\f$ satisfies the
$\Delta_2$-condition at $0$ iff the finite sequences are dense in
$\el_\f$ (see \cite{LT-I}).

\begin{thm}\label{c1s orlicz}
 Let $E, F\in HQS(C\op R)$ be infinite dimensional.
Then $\Pi_1^o(E,\, F)\subset S_2$ and there exists  an Orlicz
function $\f$ satisfying the $\Delta_2$-condition such that
$\Pi_1^o(E,\, F)=S_\f$ isomorphically. Moreover, all relevant
constants depend only on the homogeneity constants of $E$ and $F$.
 \end{thm}

We require the following lemma for the proof. Let us recall that
the noncommutative $L_1(\M)$ associated with a von Neumann algebra
$\M$ is equipped with its natural operator space structure.
Namely, as operator space, $L_1(\M)$ is the predual of the
opposite von Neumann algebra $\M^{\rm op}$ of $\M$ (see
\cite[Chapiter~7]{pis-intro}). Let $E$ and $F$ be two subspaces of
$L_1(\M)$. We denote by $E\ot_1F$ the closure of $E\ot F$ in
$L_1(\M\bar\ot\M)$. As usual, we identify an element of $E\ot F$
with its associated map from $E^*$ to $F$. Recall that a von
Neumann algebra is called QWEP if it is a quotient of a von
Neumann algebra with the WEP (weak expectation property).

\begin{lem}\label{pi1-L1}
 Let $\M$ be a QWEP von Neumann algebra and let $E, F\subset
L_1(\M)$ be finite dimensional subspaces. Then $\Pi_1^o(E^*,\,
F)=E\ot_1 F$ isometrically.
 \end{lem}

This is \cite[Corollary~10]{yew}. It is also implicitly contained
in \cite{ju-OH}. Note that we will need this lemma only for
injective $\M$. In this latter case, it immediately follows from
\cite[Lemma~5.14]{pis-ast}.

\medskip

\n{\it Proof of Theorem~\ref{c1s orlicz}.}  By Corollary \ref{hom
3 spaces}, we have
 $$C\cap R\subset E,\, F\subset C+R.$$
Consequently, $\Pi_1^o(E,\, F)\subset \Pi_1^o(C\cap R,\, C+R)$.
However, it is easy to check that $\Pi_1^o(C\cap R,\, C+R)$
coincides with $S_2$ isometrically (see case~4 below). Thus
$\Pi_1^o(E,\, F)\subset S_2$. Therefore, by Remark \ref{c1s
ideal}, $\Pi_1^o(E,\, F)=S_\Phi$ for some symmetric sequence space
$\Phi$. The main part of the proof is, of course, to show that
$\Phi$ is an Orlicz sequence space.

By Corollary \ref{hom 3 spaces} once more, $E^*$ and $F$ are
completely isomorphic to $C$, $R$, $C\cap R$ or spaces of the form
$K_{u_c,\,u_r}$ for $E^*$ belongs to $HQS(C\op R)$ too. We will
consider different cases according to the isomorphism types of
$E^*$ and $F$. The most important one is where both $E^*$ and $F$
are of type $K_{u_c,\,u_r}$.

\medskip\n{\it Case~1: $E^*=K_{u_c,\,u_r}$ and $F=K_{v_c,\,v_r}$ for
two pairs of weights on $(\O,\nu)$.} We will denote by $(e_k)$ the
canonical basis of both $K_{u_c,\,u_r}$ and $K_{v_c,\,v_r}$ (so
$E^*$ and $F$ are identified with $\el_2$ at the Banach space
level). By \cite{xu-embed}, $G_{u_c,\,u_r}$ and $G_{v_c,\,v_r}$
completely embed into a noncommutative $L_1(\M)$ with $\M$ a QWEP
von Neumann algebra with universal constants. We point out that
$\M$ can be further chosen to be injective by virtue  of
\cite{ju-araki} (see also \cite{haag-mu-kh} for a different
approach). Without loss of generality assume that all these spaces
are subspaces of $L_1(\M)$.  Then using Lemma \ref{pi1-L1} we
deduce
 $$\Pi_1^o(E,\, F)=E^* \ot_1 F\subset
 G_{u_c,\,u_r}\ot_1 G_{v_c,\,v_r}\,.$$
Indeed, since $E$ and $F$ are homogeneous and Hilbertian, we can
easily reduce the equality above to the finite dimensional case.
Assume further that all column and row spaces considered below are
subspaces of $L_1(\M)$. By the definition of $G_{u_c,\,u_r}$, we
have
 \be
 G_{u_c,\,u_r}\ot_1 G_{v_c,\,v_r}
 &=&L_2^c(u_c;\el_2)\ot_1L_2^c(v_c;\el_2)
 +L_2^r(u_r;\el_2)\ot_1L_2^r(v_r;\el_2)+\\
 &&L_2^c(u_c;\el_2)\ot_1L_2^r(v_r;\el_2)
 +L_2^r(u_r;\el_2)\ot_1L_2^c(v_c;\el_2)\,.
 \ee
Note that for two Hilbert spaces $H$ and $K$
 $$H^c\ot_1 K^c= (H\ot_2 K)^c\quad\mbox{and}\quad
 H^c\ot_1K^r=H^c\wh\ot K^r=S_1(H^*,\, K)$$
hold completely isometrically, where $\ot_2$ denotes the Hilbert
space tensor product. Recall that $S_1(H^*,\, K)$ is the trace
class of operators from $H^*$ to $K$.  Thus
 $$L_2^c(u_c;\el_2)\ot_1L_2^c(v_c;\el_2)
 =L_2^c(u_c\ot v_c;\el_2(\nz^2))
 =L_2^c(u_c\ot v_c;S_2)\,,$$
where $u_c\ot v_c$ denotes the tensor product weight of $u_c$ and
$v_c$ on $(\O\times \O,\,\nu\ot\nu)$. Similarly,
 $$L_2^r(u_r;\el_2)\ot_1L_2^r(v_r;\el_2)
 =L_2^r(u_r\ot v_r;S_2)\,.$$
On the other hand,
 $$L_2^c(u_c;\el_2)\ot_1L_2^r(v_r;\el_2)
 =L_2^c(u_c)\wh\ot L_2^r(v_r)\wh\ot S_1$$
and a similar formula for
 $L_2^r(u_r;\el_2)\ot_1L_2^c(v_c;\el_2)\,.$
Combining the preceding equalities we deduce
 \begin{eqnarray}\label{G formula}
 \begin{array}{ccl}
 \begin{displaystyle}
 G_{u_c,\,u_r}\ot_1 G_{v_c,\,v_r}\end{displaystyle}
 &=&\begin{displaystyle}L_2^c(u_c\ot v_c;S_2)+
 L_2^r(u_r\ot v_r;S_2) +\end{displaystyle}\\
 &&\begin{displaystyle}\left(L_2^c(u_c)\wh\ot L_2^r(v_r) +
 L_2^r(u_r)\wh\ot L_2^c(v_c)\right)\widehat\ot S_1\end{displaystyle}
 \end{array}
 \end{eqnarray}
completely isomorphically. We now consider these spaces only at
the Banach space level. It is easy to check that
 $$
 L_2(u_c\ot v_c;X)+
  L_2(u_r\ot v_r;X)=
  L_2(\min(u_c\ot v_c,\;u_r\ot v_r); X)
 $$
holds for any Banach space $X$ with universal equivalence
constants. Consequently,
 \beq\label{S2 term}
 L_2^c(u_c\ot v_c;S_2)+
 L_2^r(u_r\ot v_r;S_2)=A\ot_2S_2\,,
  \eeq
where
 $$A=L_2(\min(u_c\ot v_c,\;u_r\ot v_r)).$$
Let $\ot_\pi$ denote the Banach space projective tensor product.
We have
 \beq\label{S1 term}
 L_2^c(u_c)\wh\ot L_2^r(v_r) +
 L_2^r(u_r)\wh\ot L_2^c(v_c)
 =L_2(u_c)\ot_\pi L_2(v_r) +
 L_2(u_r)\ot_\pi L_2(v_c) \,{\mathop =^{\rm def}}\; B.
 \eeq

Now let  $x\in\Pi_1^o(E,\, F)$. By Remark \ref{c1s ideal}, we can
assume that $x$  is a diagonal operator:
 $$x=\sum_kx_ke_k\ot e_k\,.$$
Then by \eqref{G formula} we deduce that
 \beq\label{c1s norm matrix}
 \pi_1^o(x)\sim_c
 \inf\big\{ \|y\|_{L_2^c(u_c\ot v_c;S_2)+
 L_2^r(u_r\ot v_r;S_2)}
 + \|z\|_{\left(L_2^c(u_c)\wh\ot L_2^r(v_r) +
 L_2^r(u_r)\wh\ot L_2^c(v_c)\right)\widehat\ot S_1}\big\}\,,
 \eeq
where the infimum runs over all decompositions
$x_k=y_{kk}(\o,\s)+z_{kk}(\o,\s)$ for almost all
$(\o,\s)\in\O\times\O$ and  all $k\in\nz$ with
 \be
 y&=&(y_{kl})\in L_2^c(u_c\ot v_c;S_2)+
 L_2^r(u_r\ot v_r;S_2)\,,\\
 z&=&(z_{kl})\in \left(L_2^c(u_c)\wh\ot L_2^r(v_r) +
 L_2^r(u_r)\wh\ot L_2^c(v_c)\right)\widehat\ot S_1\,.
 \ee
Here we have viewed the elements  in $L_2^c(u_c\ot v_c;S_2)+
 L_2^r(u_r\ot v_r;S_2)$ as matrices with entries in
$L_2^c(u_c\ot v_c)+ L_2^r(u_r\ot v_r)$ and similarly for the other
space. Note that the diagonal projection on $S_2$ extends to
contractive projections  on both
 $L_2^c(u_c\ot v_c;S_2)+ L_2^r(u_r\ot v_r;S_2)$
and
 $\left(L_2^c(u_c)\wh\ot L_2^r(v_r) +
 L_2^r(u_r)\wh\ot L_2^c(v_c)\right)\widehat\ot S_1\,.$
Therefore, the infimum in \eqref{c1s norm matrix} can be
restricted to all decompositions with diagonal $y$ and $z$.
However, for a diagonal $y$  using \eqref{S2 term} we find
 \be
  \|y\|_{L_2^c(u_c\ot v_c;S_2)+
 L_2^r(u_r\ot v_r;S_2)}
  \sim\|y\|_{\el_2(A))}
  =\big(\sum_k \|y_k\|^2_{A}\big)^{1/2}\,.
  \ee
On the other hand, if $z$ is diagonal, we have, by virtue of
\eqref{S1 term}
 \be
 \|z\|_{\left(L_2^c(u_c)\wh\ot L_2^r(v_r) +
 L_2^r(u_r)\wh\ot L_2^c(v_c)\right)\widehat\ot S_1}
 &=& \|z\|_{\left(L_2^c(u_c)\wh\ot L_2^r(v_r) +
 L_2^r(u_r)\wh\ot L_2^c(v_c)\right)\widehat\ot\el_1}\\
 &=&\|z\|_{\el_1(L_2^c(u_c)\wh\ot L_2^r(v_r) +
 L_2^r(u_r)\wh\ot L_2^c(v_c))}\\
 &=&\|z\|_{\el_1(B)}\,.
 \ee
It thus follows that
 \beq\label{c1s norm diag}
 \pi_1^o(x)\sim_c
 \inf\big\{\big(\sum_k \|y_k\|^2_{A}\big)^{1/2}
 + \sum_k \|z_k\|_{B}\big\}\,,
 \eeq
where the infimum runs over all decompositions
$x_k=y_k(\o,\s)+z_k(\o,\s)$ for almost all $(\o,\s)\in\O\times\O$
and  all $k\in\nz$ with $y_k\in A$ and $z_k\in B$. This is the key
formula of the proof. It suggests the definition of the desired
Orlicz function $\f$. Indeed, define $\wt\f: \rz_+\to\rz_+$ by
 \be
 \wt\f(t)=\inf\big\{t^2\|a\|^2_{A}+ t\,\|b\|_{B}\,:\,
 a\in A,\, b\in B\;{\rm s.t.}\;a+b=1\;{\rm a.e.\; on}\;
 \O\times\O\big\}.
 \ee
Note that $\wt\f$ satisfies the following properties:
 \begin{enumerate}[(i)]
 \item $\wt\f$ is continuous, $\wt\f(0)=0$ and $\lim_{t\to\8}\wt\f(t)=\8$;
 \item both $\wt\f$ and $\wt\f/t$ are nondecreasing;
 \item $\wt\f(2t)\le 4\,\wt\f(t)$, so $\wt\f$ verifies
the $\Delta_2$-condition.
 \end{enumerate}
Let
 $$\f(t)=\int_0^t\frac{\wt\f(s)}s\,ds.$$
Then $\f$ is an Orlicz function satisfying the
$\Delta_2$-condition. Moreover,
 $$\f(t)\le\wt\f(t)\le 4\,\f(t),\quad\forall\;t\in\rz_+.$$
Therefore $\f$ is equivalent to $\wt\f$.

It is now easy to show that $\pi_1^o(x)\sim_c\|x\|_\f$ for any
$x=(x_k)$. Assume $\|x\|_\f<1$. Then there exist $(a_k)$ and
$(b_k)$ such that $a_k+b_k=1$ a.e. on $\O\times\O$ and
 $$\sum_k \big(|x_k|^2\,\|a_k\|_A^2 + |x_k|\,\|b_k\|_B\big)\le 4.$$
Set $y_k=x_ka_k$ and $z_k=x_kb_k$. Then $x_k=y_k+z_k$ a.e. on
$\O\times\O$ and
 $$\sum_k \|y_k\|^2_{A}+\sum_k \|z_k\|_{B}\le 4\,.$$
Thus by \eqref{c1s norm diag}, $\pi_1^o(x)\le c$. Conversely, let
$x$ be such that the infimum in \eqref{c1s norm diag} is less than
$1$. Choose two sequences $(y_k)$ and $(z_k)$ such that
$x_k=y_k+z_k$ and
 $$\big(\sum_k \|y_k\|^2_{A}\big)^{1/2} + \sum_k \|z_k\|_{B}<1.$$
Then
 $$\sum_k \big(|x_k|^2\,\|a_k\|_A^2 + |u_k|\,\|b_k\|_B\big)<1,$$
where $a_k=y_k/x_k$ and $b_k=z_k/x_k$. Since $a_k+b_k=1$, we find
 $$\sum_k\f(|x_k|)\le \sum_k\wt\f(|x_k|)<1;$$
whence $\|x\|_\f\le 1$. Therefore, we get the desired equivalence
$\pi_1^o(x)\sim_c\|x\|_\f$. This shows the theorem in the case
where $E^*=K_{u_c,\,u_r}$ and $F=K_{v_c,\,v_r}$.

\medskip\n{\it Case~2: $E^*\in\{C,\,R,\,C\cap R\}$
and $F=K_{v_c,\,v_r}$.} The proof in this case is similar to the
previous one but much simpler. Consider first the case where
$E^*=C$. Using the identification $C\simeq R_1$, we have
 \be
 C\ot_1 G_{v_c,\,v_r}=R_1[ G_{v_c,\,v_r}]
 =R_1[L_2^c(v_c;\el_2)]+ R_1[L_2^r(v_r;\el_2)]
 =L_2^c(v_c;S_2) + L_2^r(v_r)\wh\ot S_1\,.
 \ee
This is the analogue of \eqref{G formula} for the present case.
Then as before, for a diagonal operator $x\in\Pi_1^o(E,\,F)$ we
have
 \be
 \pi_1^o(x)\sim_c
 \inf\big\{\big(\sum_k \|y_k\|^2_{L_2(v_c)}\big)^{1/2}
 + \sum_k \|z_k\|_{L_2(v_r)}\big\}\,,
 \ee
where the infimum runs over all decompositions
$x_k=y_k(\o)+z_k(\o)$ for almost all $\o\in\O$ and $k\in\nz$ with
$(y_k)\subset L_2(v_c)$ and $(z_k)\subset L_2(v_r)$. Like in
case~1, we deduce that $\pi_1^o(x)\sim_c\|x\|_\f$, where $\f$ is
an Orlicz function equivalent to
 $$\wt\f(t)=\inf\big\{t^2\|a\|^2_{L_2(v_c)}+
 t\,\|b\|_{L_2(v_r)}\,:\,
 a\in L_2(v_c),\, b\in L_2(v_r)\;{\rm s.t.}\;a +b=1\;
 {\rm a.e.\; on}\; \O\big\}.$$
The case $E^*=R$ is dealt with similarly. Now assume $E^*=C\cap
R$. Then
 \be
 (C\cap R)\ot_1 G_{v_c,\,v_r}=
 \big(L_2^c(v_c;S_2) + L_2^r(v_r)\wh\ot S_1\big)\cap
 \big(L_2^c(v_c)\wh\ot S_1 + L_2^r(v_r;S_2)\big)\,.
 \ee
Therefore for a diagonal operator $x\in\Pi_1^o(E,\,F)$ we deduce
 \be
 \pi_1^o(x)\sim_c
 \inf\Big\{\big[\sum_k \big(\|y_k\|^2_{L_2(v_c)}
 +\|z'_k\|^2_{L_2(v_r)}\big)\big]^{1/2}
 + \sum_k \big(\|y_k'\|_{L_2(v_c)}+\|z_k\|_{L_2(v_r)}\big) \Big\}\,,
 \ee
where the infimum runs over all decompositions
$x_k=y_k+z_k=y'_k+z'_k$ a.e. on $\O$ for $k\in\nz$ with
$(y_k),\,(y'_k)\subset L_2(v_c)$ and $(z_k),\, (z'_k)\subset
L_2(v_r)$. This time, the required Orlicz function $\f$ is
equivalent to the following one
 \be
 &&\wt\f(t)=\inf\big\{t^2\big(\|a\|^2_{L_2(v_c)}
 +\|b'\|^2_{L_2(v_r)}\big)+
 t\big(\|a'\|_{L_2(v_c)}+\|b\|_{L_2(v_r)}\big)\,:\,\\
 &&\hskip 1.8cm a,\,a'\in L_2(v_c),\, b,\,b'\in L_2(v_r)\;{\rm s.t.}\;
 a+b=1=a'+b'\big\}.
 \ee

\medskip\n{\it Case~3: $E^*=K_{u_c,\,u_r}$ and
$F\in\{C,\,R,\, C\cap R\}$.} This
is symmetric to case~2 and is proved in the same way. We thus omit
the details.

\medskip\n{\it Case~4: $E^*,\,F\in\{C,\,R,\, C\cap R\}$.}
This is the trivial case since we have the following easily
checked isometric identifications
 \be
 &&\Pi_1^o(C,\, C)=\Pi_1^o(R,\, R)=S_1,\quad
 \Pi_1^o(C,\, R)=\Pi_1^o(C,\, R)=S_2\,,\\
  &&\Pi_1^o(C+R,\, C)=\Pi_1^o(C+R,\, R)=\Pi_1^o(C+R,\, C\cap R)=S_1\,.
 \ee
Let us check, for instance, $\Pi_1^o(C+R,\, C)=S_2$. We have
 $$\Pi_1^o(C+R,\, C)=(C\cap R)\ot_1 C=(R_1\cap C_1)\ot_1 R_1
 =R_1[R_1]\cap C_1[R_1]= S_2\cap S_1= S_1\,.$$
Therefore, the proof of the theorem is complete. \cqd

\medskip

Since an Orlicz function $\f$ is uniquely determined by its
fundamental sequence $(\f_n)_{n\ge1}$
 $$\f_n=\big\|\sum_{k\le n}e_k\big\|_\f=\frac1{\f^{-1}(\frac1n)}\,,$$
Theorem \ref{c1s orlicz} reduces the calculation of the space
$\Pi_1^o(E,\, F)$ to that of $\pi_1^o(\id_n: E_n\to F_n)$ for all
$n$, where $E_n$ and $F_n$ are the linear spans of $\{e_1, ...,
e_n\}$ in $E$ and $F$, respectively, and where $\id_n: E_n\to F_n$
denotes the formal identity from $E_n$ to $F_n$, i.e.,
$\id_n(e_k)=e_k$ for every $k=1, ..., n$. We record this in the
following

\begin{cor}\label{c1s funda}
 With the notations above  $\Pi_1^o(E,\,
F)$ is uniquely determined, up to isomorphism, by the sequence
$(\pi_1^o(\id_n\,:\, E_n\to F_n))_{n\ge1}$.
 \end{cor}

This result considerably simplifies the task of determining
$\Pi_1^o(E,\, F)$. Thanks to Theorem~\ref{c1s orlicz}, the only
things left to be calculated are the completely 1-summing norms of
the identities. For these we have a more concrete formula:

\begin{prop}\label{c1s id}
 Let $(u_c,\,u_r)$ and $(v_c,\, v_r)$ be two pairs of weights on
$(\O,\nu)$ satisfying \eqref{w-con}. Let $E^*=K_{u_c,\,u_r}$ and
$F=K_{v_c,\,v_r}$. Then
 \beq\label{c1s norm id}
 \pi_1^o(\id_n)\sim_c
 \inf\left\{\sqrt n\,\|a\|_{A}
 + n\,\|b\|_{B}\,:\, a\in A,\; b\in B,\; a+b=1\;
 {\rm a.e.\; on}\; \O\right\}.
 \eeq
where
 \be
 A=L_2(\min(u_c\ot v_c,\;u_r\ot v_r)), \quad B=L_2(u_c)\ot_\pi L_2(v_r)+
 L_2(u_r)\ot_\pi L_2(v_c).
 \ee

  \end{prop}

\pf Consider the projection $P_n$ defined for sequences $x=(x_k)$
by
 $$P_n(x)=\frac{x_1+\cdots +x_n}n\,\sum_{k=1}^ne_k\,.$$
Then $P_n$ extends to a contractive projection on both $\el_2(A)$
and $\el_1(B)$. Together with \eqref{c1s norm diag}, this implies
\eqref{c1s norm id}.
 \cqd

\medskip

The remainder of this section is essentially devoted to the proof of
Theorem~\ref{c1s id explicit}.

\medskip\n\emph{Proof of Theorem~\ref{c1s id explicit}.}
The first part of the theorem is already contained in Theorem
\ref{c1s orlicz}. We will prove the estimate on the fundamental
sequence of $\f$. Recall that the fundamental functions $\F_{c,
E^*}, \F_{r, E^*}, \F_{c, F}$ and $\F_{r, F}$ are extended to
continuous functions on $\rz_+$. Without loss of generality wa can
further assume that they are increasing on $\rz_+$ and equal $1$
at $0$. Let
  $$
  w_{c, F}(t)=\left\{\begin{array}{ll}
 \displaystyle 1 & \textrm{ if } t>0\\
 \displaystyle\frac{1}{\F_{c, F}^{-1}(-t)}&\textrm{ if } t\le0
 \end{array}\right.\quad\mbox{and}\quad
 w_{r, F}(t)=\left\{\begin{array}{ll}
 \displaystyle 1 & \textrm{ if } t<0\\
 \displaystyle\frac{1}{\F_{r, F}^{-1}(t)}&\textrm{ if } t\ge0
 \end{array}\,.\right.
 $$
The pair $(w_{c, E^*},\, w_{r, E^*})$ has the same meaning with
$E^*$ in place of $F$.  Then by Theorem \ref{c1s orlicz}
 $$E^*=K_{w_{c, E^*},\, w_{r, E^*}}\quad\mbox{and}\quad
 F=K_{w_{c, F},\, w_{r, F}}\,.$$
Now we are in a position to apply \eqref{c1s norm id}. According
to the four quadrants of the plane, we decompose the space $A$
there as follows
 \be
 && A=
 L_2(\min(w_{c, E^*}^-\ot w_{c, F}^-\,,\;
 w_{r, E^*}^-\ot w_{r, F}^-))
 \op L_2(\min(w_{c, E^*}^+\ot w_{c, F}^+\,,
 \;w_{r, E^*}^+\ot w_{r, F}^+))\\
 &&\,\hskip .3cm \op\; L_2(\min(w_{c, E^*}^-\ot w_{c, F}^+\,,
 \;w_{r, E^*}^-\ot w_{r, F}^+))
 \op L_2(\min(w_{c, E^*}^+\ot w_{c, F}^-\,,
 \;w_{r, E^*}^+\ot w_{r, F}^-))\\
 &&\hskip .3cm{\displaystyle\mathop=^{\rm def}}\; A_{--} \op A_{++}
 \op A_{-+}\op A_{+-}\,,
 \ee
where
 $$w^-=w\big|_{\rz_-} \quad\mbox{and}\quad
  w^+=w\big|_{\rz_+}\,.$$
A similar decomposition holds for $B$. Consequently, the infimum
in \eqref{c1s norm id} is a sum of four infima corresponding to
these decompositions. More precisely, we have
 \be
 \pi_1^o(\id_n)\sim
 \|1\|_{\sqrt n\,A_{--}\, +\,n B_{--}} +
 \|1\|_{\sqrt n\,A_{++}\, +\, n B_{++}}
 + \|1\|_{\sqrt n\,A_{-+}\, +\,n B_{-+}}
 + \|1\|_{\sqrt n\,A_{+-}\, +\,n B_{+-}}\,.
 \ee
Here and during this proof all equivalence constants depend only
on the homogeneity and regularity constants of $E$ and $F$. Thus
it suffices to deal with  separately the four terms on the right.
The first two are easy. Indeed, we have
 \be
 A_{--}
 =L_2(\min(w_{c, E^*}^-\ot w_{c, F}^-\,,\;
 1\ot 1))
 =L_2(w_{c, E^*}^-\ot w_{c, F}^-).
 \ee
Since $w_{c, E^*}^-\ot w_{c, F}^-$ is integrable on $\rz_-\times
\rz_-$, $1\in L_2(w_{c, E^*}^-\ot w_{c, F}^-)$. Thus
 $$\|1\|_{\sqrt n\,A_{--}\, +\,n B_{--}}\sim\sqrt n\,.$$
Similarly,
 $$\|1\|_{\sqrt n\,A_{++}\, +\,n B_{++}}\sim\sqrt n\,.$$
Therefore,
 \beq\label{-++-}
 \pi_1^o(\id_n)\sim
  \|1\|_{\sqrt n\,A_{-+}\, +\,n B_{-+}}
 + \|1\|_{\sqrt n\,A_{+-}\, +\,n B_{+-}}\,.
 \eeq
To estimate the last two terms, note that
 $$ A_{-+}=  L_2(\min(w_{c, E^*}^-\ot 1\,,
 \;1\ot w_{r, F}^+))$$
and
 $$B_{-+}=L_2(w_{c, E^*}^-)\ot_\pi  L_2(w_{r, F}^+)+
 L_2(\rz_-)\ot_\pi L_2(\rz_+)
 =L_2(w_{c, E^*}^-)\ot_\pi  L_2(w_{r, F}^+).$$
Since the projective tensor norm dominates the Hilbert tensor
norm, we have the following contractive inclusion
 $$B_{-+}\subset L_2(w_{c, E^*}^-)\ot_2  L_2(w_{r, F}^+)
 =L_2(w_{c, E^*}^-\ot w_{r, F}^+).$$
Hence,
 $$\|1\|_{\sqrt n\,A_{-+}\, +\,n B_{-+}}\ge c^{-1}
 \|1\|_{L_2(\min(n\,w_{c, E^*}^-\ot 1\,,
 \;n\,1\ot w_{r, F}^+\,,\; n^2\,w_{c, E^*}^-\ot w_{r, F}^+))}\,.$$
To ease the calculation of the minimum above we set $\wt w_{c,
E^*}(s)=w_{c, E^*}^-(-s)$ for $s\in\rz_+$. Since $\wt w_{c, E^*}$
and $w_{r, F}^+$ are decreasing on $\rz_+$,  for each $s\in\rz_+$
there exists a unique $t=t(s)\in\rz_+$ such that
 $$\wt w_{c, E^*}(s)=w_{r, F}^+(t(s)).$$
The function $s\mapsto t(s)$ is  increasing and bijective on
$\rz_+$. Its inverse is the function $t\mapsto s(t)$ such that
 $$\wt w_{c, E^*}(s(t))=w_{r, F}^+(t).$$
We then have
 $$\min(\wt w_{c, E^*}(s),\;w_{r, F}^+(t))=
 \left\{\begin{array}{ll}
 \displaystyle \wt w_{c, E^*}(s) & \textrm { if } t\le t(s),\\
 \displaystyle w_{r, F}^+(t) & \textrm{ if } t> t(s).
 \end{array}\right.$$
On the other hand, for each $n\in\nz$ there exist (unique) $s_n,
t_n\in\rz_+$ such that
 $$\wt w_{c, E^*}(s_n)=w_{r,F}^+(t_n)=\frac1n\,.$$
Note that $t_n=t(s_n)$. We now decompose $\rz_+\times \rz_+$ into
a union of three disjoint regions:
 \be
 \La_1&=&\left\{(s, t)\,:\,  t\le t_n,\; s\ge s(t)\right\},\\
 \La_2&=&\left\{(s, t)\,:\, s\le s_n,\; t>t(s)\right\},\\
 \La_3&=&(s_n,\;\8)\times (t_n,\;\8).
 \ee
Then
 $$\min\big(n\wt w_{c, E^*}(s),\;nw_{r, F}^+(t)\,,\;
 n^2\wt w_{c, E^*}(s)w_{r, F}^+(t)\big)=
 \left\{\begin{array}{ll}
 \displaystyle n\wt w_{c, E^*}(s) & \textrm { if } (s, t)\in\La_1\,,\\
 \displaystyle nw_{r, F}^+(t) & \textrm{ if } (s, t)\in\La_2\,,\\
 \displaystyle n^2\,\wt w_{c, E^*}(s) w_{r, F}^+(t) & \textrm{ if }
  (s, t)\in\La_3\,.
 \end{array}\right.$$
Therefore,
 \be
 &&\|1\|^2_{L_2(\min(n\,w_{c, E^*}^-\ot 1\,,
 \;n\,1\ot w_{r, F}^+\,,\; n^2\,w_{c, E^*}^-\ot w_{r, F}^+))}\\
 &&~~=n\int_{\La_1}\wt w_{c, E^*}(s) dsdt
 +n\int_{\La_2}w_{r, F}^+(t) dsdt
 + n^2\int_{\La_3}\wt w_{c, E^*}(s) w_{r, F}^+(t) dsdt\\
 &&~~{\displaystyle\mathop =^{\rm def}}\;\l_1 +\l_2+\l_3.
 \ee
Recall that $h_{c, E^*}$ and $g_{c, E^*}$ are the functions
associated to $\wt w_{c, E^*}$ by \eqref{def hg}. Then by
\eqref{g-h}
 \be
 \frac{\l_3}{n^2}
 &=&\int_{s_n}^\8\wt w_{c, E^*}(s)ds\,
 \int_{t_n}^\8 w_{r, F}^+(t) dt\\
 &=&h_{c, E^*}(s_n)h_{r, F}(t_n)
 =\frac{s_nt_n}{g^{-1}_{c, E^*}(s_n)\, g^{-1}_{r,
 F}(t_n)}\,.
 \ee
However, by \eqref{g-w}
 \be
 \frac{1}{g^{-1}_{c, E^*}(s_n)}\sim \wt w_{c,
 E^*}(s_n)=\frac1n\,;
 \ee
whence
 $$s_n\sim g_{c, E^*}(n).$$
Similar estimates hold for $g^{-1}_{r,
 F}(t_n)$. It follows that
 $$\l_3
 \sim g_{c, E^*}(n)g_{r,F}(n).$$
For $\l_2$ we find
 \be
 \frac{\l_2}n
 =\int_0^{s_n}\int_{t(s)}^\8w_{r, F}^+(t)\, dtds
 =\int_0^{s_n}h_{r, F}(t(s))ds
 =\int_0^{s_n}\frac{t(s)}{g^{-1}_{r, F}(t(s))}ds\,.
 \ee
However,
 $$\frac1{g^{-1}_{r, F}(t(s))}\sim
 w_{r, F}^+(t(s))=\wt w_{c, E^*}(s)
 \sim \frac1{g^{-1}_{c, E^*}(s)}\,.$$
Thus
 $$t(s)\sim g_{r, F}(g^{-1}_{c, E^*}(s)).$$
Consequently,
 $$
 \l_2
 \sim n \int_0^{g_{c, E^*}(n)}
 \frac{g_{r, F}(g^{-1}_{c, E^*}(s))}{g^{-1}_{c, E^*}(s)}\,ds.
 $$
The same calculation applies to $\l_1$ and yields
 $$
 \l_1
 \sim n \int_0^{g_{r, F}(n)}
 \frac{g_{c, E^*}(g^{-1}_{r, F}(t))}{g^{-1}_{r, F}(t)}\,dt.
 $$
Combining the preceding estimates we obtain
 \be
 && c^2
 \|1\|^2_{\sqrt n\,A_{-+}\, +\,n B_{-+}}\ge g_{c, E^*}(n)g_{r,F}(n)+\\
 &&~~ n\int_0^{g_{c, E^*}(n)}
 \frac{g_{r, F}(g^{-1}_{c, E^*}(s))}{g^{-1}_{c, E^*}(s)}\,ds+
 n \int_0^{g_{r, F}(n)}
 \frac{g_{c, E^*}(g^{-1}_{r, F}(t))}{g^{-1}_{r, F}(t)}\,dt\,.
 \ee
To prove the converse inequality it suffices to note  that the
right hand side above corresponds to the decomposition $1=a+b$
with
 $$a=\un_{\La_1} + \un_{\La_2}\in A_{-+}\quad\mbox{and}\quad
 b=\un_{\La_3}\in B_{-+}\,.$$
We have
 $$\sqrt n\,\|a\|_{A_{-+}}=\sqrt{\l_1+\l_2}\quad\mbox{and}\quad
 n\|b\|_{ B_{-+}}=\sqrt{\l_3}\,.$$
Although not needed, it is helpful to observe that $b$ is a
tensor, so its norm in $B_{-+}$ (the projective norm) coincides
with its Hilbert tensor norm. Thus
 \be
 \|1\|_{\sqrt n\,A_{-+}\, +\,n B_{-+}}\le
 \sqrt{\l_1+\l_2}+ \sqrt{\l_3}\,.
 \ee
Therefore, we deduce
 \begin{eqnarray}\label{-+}
 \begin{array}{ll}
  &\displaystyle\|1\|^2_{\sqrt n\,A_{-+}\, +\,n B_{-+}}\sim
 g_{c, E^*}(n)g_{r,F}(n)\;+\\
 &~ \displaystyle n\int_0^{g_{c, E^*}(n)}
 \frac{g_{r, F}(g^{-1}_{c, E^*}(s))}{g^{-1}_{c, E^*}(s)}\,ds+
 n \int_0^{g_{r, F}(n)}
 \frac{g_{c, E^*}(g^{-1}_{r, F}(t))}{g^{-1}_{r, F}(t)}\,dt\,.
 \end{array}
 \end{eqnarray}
By symmetry, we also find
 \begin{eqnarray}\label{+-}
 \begin{array}{ll}
  &\displaystyle\|1\|^2_{\sqrt n\,A_{+-}\, +\,n B_{+-}}\sim
 g_{r, E^*}(n)g_{c,F}(n)\;+\\
 &~ \displaystyle n\int_0^{g_{r, E^*}(n)}
 \frac{g_{c, F}(g^{-1}_{r, E^*}(s))}{g^{-1}_{r, E^*}(s)}\,ds+
 n \int_0^{g_{c, F}(n)}
 \frac{g_{r, E^*}(g^{-1}_{c, F}(t))}{g^{-1}_{c, F}(t)}\,dt\,.
 \end{array}
 \end{eqnarray}
Combining \eqref{-++-}, \eqref{-+} and \eqref{+-}, we finally get
 \be
 &&\pi_1^o(\id_n)^2
 \sim g_{c, E^*}(n)g_{r,F}(n) + g_{r, E^*}(n)g_{c,F}(n) +\\
 &&\hskip 1.6cm n\int_0^{g_{c, E^*}(n)}
 \frac{g_{r, F}(g^{-1}_{c, E^*}(t))}{g^{-1}_{c, E^*}(t)}\,dt+
 n \int_0^{g_{r, F}(n)}
 \frac{g_{c, E^*}(g^{-1}_{r, F}(t))}{g^{-1}_{r, F}(t)}\,dt +\\
 &&\hskip 1.6cm n\int_0^{g_{r, E^*}(n)}
 \frac{g_{c, F}(g^{-1}_{r, E^*}(t))}{g^{-1}_{r, E^*}(t)}\,dt+
 n \int_0^{g_{c, F}(n)}
 \frac{g_{r, E^*}(g^{-1}_{c, F}(t))}{g^{-1}_{c, F}(t)}\,dt\,.
 \ee
By the definition of the weights in consideration and \eqref{g-w}
we have
 $g_{c, E^*}\sim \F_{c, E^*}$
and similar equivalences for all other functions. This allows us
to replace the functions $g$ above by the respective fundamental
functions. On the other hand, since all functions in consideration
are equivalent to $1$ on $(0, 1)$, the parts on $(0,1)$ of all
integrals above can be disregarded. Finally, recalling
$\f_n\sim\pi_1^o(\id_n)$, we complete the proof of the theorem.
\cqd

\begin{rk}
 Assume  that the two pairs of weights in
Proposition~\ref{c1s id} satisfy \eqref{reg1}. Then the preceding
proof shows that the infimum in \eqref{c1s norm id} can be
restricted to indicator functions $a$ and $b$.
 \end{rk}

Theorem~\ref{c1s id explicit} clearly implies that $\pi_1^o(x)\sim
\pi_1^o(x^*)$ for any $x\in\Pi_1^o(E,\, F)$ and $E, F\in HQS(C\op
R)$. This fact is true for general $E, F\in QS(C\op R)$:

\begin{prop}\label{c1s-sym}
 Let $E, F\in QS(C\op R)$. Then $\pi_1^o(x)\sim \pi_1^o(x^*)$
for any $x\in\Pi_1^o(E,\, F)$ with universal equivalence
constants.
\end{prop}

\pf Since $E$ and $F$ have the completely bounded approximation
property with a universal constant, we can assume that both $E$
and $F$ are finite dimensional. On the other hand, recall that any
space $QS(C\op R)$ completely embeds into a noncommutative
$L_1(\M)$ with $\M$ a QWEP. It then remains to apply
Lemma~\ref{pi1-L1}. \cqd


\section{Injectivity and exactness}
 \label{Injectivity and exactness}


Recall that an operator space $F$ is called {\it injective} if the
identity map of $F$ factors through $B(H)$ by completely bounded
maps for some Hilbert space $H$, or equivalently, if $F$ is
completely complemented in $B(H)$. Let $F$ be a (completely
isometric) subspace of $B(H)$ for some Hilbert space $H$. The {\it
projection} or {\it injectivity constant} of $F$ is then defined
to be
  $$\l_{cb}(F) = \inf\left\{\|P\|_{cb}\,:\, F\subset B(H)
 \textrm{ as subspace, } P:B(H)\to F \textrm{
 projection}\right\}\,.$$

Using Proposition \ref{c1s-sym} and the trace duality between the
completely 1-summing norm and the $\gamma_\8$-norm (see
\cite[Lemma~4.6]{ju-OH}), we immediately deduce the following

\begin{prop}\label{inj-sym}
 Let $F\in QS(C\op R)$. Then $\l_{cb}(F)\sim \l_{cb}(F^*)$
with universal equivalence constants.
 \end{prop}

\begin{lem}\label{pi1-gama8}
 Let $F$ be an $n$-dimensional $\l$-homogeneous $\mu$-Hilbertian
operator space. Then
 $$n\le\pi_1^o(\id_F)\,\l_{cb}(F)\le\l\mu\,n.$$
 \end{lem}

\pf This lemma is proved in \cite{ju-OH} for the $n$-dimensional
operator Hilbert space $OH_n$ (see the proof of
\cite[Corollary~4.11]{ju-OH}). The proof there uses only the
$1$-homogeneity of $OH_n$; so it remains valid for general
homogeneous Hilbertian operator spaces. We omit the details.\cqd

\medskip

\n\emph{Proof of Theorem~\ref{inj}}. This is immediate from
Lemma~\ref{pi1-gama8}, Theorems~\ref{c1s id explicit} and
\ref{dual}. \cqd

\medskip

Now we turn to the exactness. Recall that the {\it exactness
constant} of an operator space $F$ is defined by
 $$ex(F)=\sup_{E\subset F,\, \dim E<\8}\,
 \inf_{G\subset \mathbb K(\ell_2)}
 d_{cb}(E,\,G).$$
$F$ is called {\it exact} if $ex(F)<\8$. We refer to
\cite{er-book}, \cite{pis-intro} and \cite{pis-ex} for more
information.

\medskip

The exactness constant of a subspace of a noncommutative $L_1$ can
be also expressed as a projection constant, as shown by the
following result. This explains why the exactness constants of
spaces in $HQS(C\op R)$ can be dealt with  similarly  as their
projection constants. Let $F$ be an operator space. Define
 $$\wt\l_{cb}(F)=\inf\|P\|_{cb}\,,$$
where the infimum runs over all von Neumann algebras  $\N$ such
that $F\subset \N$ as a completely isometric subspace and all
completely bounded projections $P: \N\to F$. Note that if we
require $\N$ to be injective in the infimum above, then we recover
the injectivity constant of $F$.

\begin{prop}\label{ex-g8}
 Let $F$ be a finite dimensional subspace of $L_1(\M)$ for
a von Neumann algebra $\M$. Then
 $$ ex(F)\le\wt\l_{cb}(F)\le c\, ex(F),$$
where $c$ is a universal positive constant.
 \end{prop}

\pf The first inequality is \cite[Corollary~17.16]{pis-intro}.
This is true without the assumption that $F\subset L_1(\M)$.  It
remains to show the second. By the operator space Grothendieck
theorem of Pisier and Shlyakhtenko \cite{pisshlyak}, the inclusion
map $\iota: F\hookrightarrow L_1(\M)$ factors through $C\op R$ by
completely bounded maps. More precisely, there exist $x\in CB(F,\,
C\op R)$ and $y\in CB(C\op R,\, L_1(\M))$ such that $\iota=yx$ and
$\|y\|_{cb}\|x\|_{cb}\le c\,ex(F)$ for some universal constant
$c$. Let $S=x(F)\subset C\op R$. Then by \cite{xu-embed}, $S$ is
completely isomorphic to a completely complemented subspace $G$ of
a von Neumann $\N$ with universal constants.  Since $F$ is
completely isomorphic to $S$ with constant $\|y\|_{cb}\|x\|_{cb}$,
we then deduce the desired inequality.\cqd

\medskip

We now pass to consider spaces in $HQS(C\op R)$. Let $(w_c, w_r)$
be a pair of weights on $\rz$ satisfying  \eqref{reg1}. Let
$F=K_{w_c,\,w_r}$. Recall that $F$ is the subspace of constant
functions of $G_{w_c,\,w_r}=L_2^c(w_c;\el_2)+L_2^r(w_r;\el_2)$.
Let $F_n$ be an $n$-dimensional subspace of $F$. By homogeneity
$F_n$ can be assumed to be the subspace of constant functions of
$G^n_{w_c,\,w_r}=L_2^c(w_c;\el_2^n)+L_2^r(w_r;\el_2^n)$. Let
$\iota: F_n\hookrightarrow G^n_{w_c,\,w_r}$ and
$q\,:\,L_2^c(w_c;\el_2^n)\op L_2^r(w_r;\el_2^n)\to
G^n_{w_c,\,w_r}$ be the natural inclusion and quotient maps,
respectively.

\begin{lem}\label{exK}
 With the notations above we have
 $$ex(F_n)\le\inf\{\|x\|_{cb}\,:\, x:F_n\to
 L_2^c(w_c;\el_2^n)\op L_2^r(w_r;\el_2^n)\;{\rm s.t. }\;
 qx=\iota\}\le c\, ex(F_n).$$

 \end{lem}

\pf Let $x:F_n\to
 L_2^c(w_c;\el_2^n)\op L_2^r(w_r;\el_2^n)$ be such that
$qx=\iota$. Then $F_n$ is completely isomorphic to  $x(F_n)$ with
constant $\|x\|_{cb}$. Since $L_2^c(w_c;\el_2^n)\op
L_2^r(w_r;\el_2^n)$ is exact with constant $1$, we deduce
$ex(F_n)\le\|x\|_{cb}$. The nontrivial part is the upper estimate.
The proof of this is similar to that of Proposition~\ref{ex-g8}.
By \cite{xu-embed}, $G^n_{w_c,\,w_r}$ is completely isomorphic to
a complemented subspace of a noncommutative $L_1$-space with
universal constants. Thus using again Pisier and Shlyakhtenko's
theorem, we see that the inclusion $\iota: F_n\hookrightarrow
G^n_{w_c,\,w_r}$ factors through $C\op R$ by completely bounded
maps. Namely, there exist $y: F_n\to C\op R$ and $z: C\op R\to
G^n_{w_c,\,w_r}$ such that $zy=\iota$ and $\|y\|_{cb}\|z\|_{cb}\le
\l\,ex(F_n)$. Passing to duals, we get a factorization
$\iota^*=y^*z^*$ of $\iota^*$ through $R\op C$. Since $R\op C$ is
injective with constant $1$, $z^*$ admits an extension $\wt z:
L_2^r(w_c^{-1};\el_2^n)\op L_2^c(w_r^{-1};\el_2^n)\to R\op C$ with
$\|\wt z\|_{cb}=\|z\|_{cb}$. Thus we have the following
commutative diagram
 $$\xymatrix{
    L_2^r(w_c^{-1};\el_2^n)\cap L_2^c(w_r^{-1};\el_2^n)
    \ar[rr]^-{\iota^*} \ar@{^{(}->}[d]  &&  F_n^* \\
    L_2^r(w_c^{-1};\el_2^n)\op L_2^c(w_r^{-1};\el_2^n) \ar[rr]_-{\wt z}
    && R\op C\ar[u]_{y^*} }$$
Dualizing this diagram and setting $x=\wt z^*y$, we get
$qx=\iota$ and $\|x\|_{cb}\le\l\,ex(F_n)$. This implies the
desired upper estimate. \cqd

\medskip

Recall that $L_2^c(w_c;\el_2^n)=\big(L_2(w_c)\ot_2\el_2^n\big)^c$
and $F_n$ is identified with $\el_2^n$ as Banach spaces. If $a\in
L_2(w_c)$, we use $a\ot\id$ to denote the map $f\mapsto a\ot f$
from $F_n$ to $L_2^c(w_c;\el_2^n)$.

\begin{lem}\label{diag}
 Let $x:F_n\to L_2^c(w_c;\el_2^n)\op L_2^r(w_r;\el_2^n)$ be such
that $qx=\iota$. Then there exist $a\in L_2(w_c)$ and $b\in
L_2(w_r)$ such that $a+b=1$ a.e. and such that
 $$q\wt x=\iota \quad\mbox{and}\quad
 \|\wt x\|_{cb}\le c\|x\|_{cb}\,,$$
where $\wt x:F_n\to L_2^c(w_c;\el_2^n)\op L_2^r(w_r;\el_2^n)$ is
defined by $\wt x=(a\ot\id\,,\, b\ot\id)$.
 \end{lem}

\pf  Let $\mathcal U$ be the unitary group of $\el_2^n$, equipped
with Haar measure. For any $u\in \mathcal U$ we consider $u$ as
maps on both $F_n$ and $L_2^c(w_c;\el_2^n)$. More precisely,
viewed as a map on
$L_2^c(w_c;\el_2^n)=\big(L_2^c(w_c)\ot_2\el_2^n\big)^c$, $u$ acts
only on the factor $\el_2^n$, so agrees with $\id\ot u$. Writing
$x=(y, z)$, we define $\wt y\,:\, F_n\to L_2^c(w_c;\el_2^n)$ by
 $$\wt y=\int_{\mathcal U}u^*yu\,du.$$
More precisely,
 $$\wt y(f)=\int_{\mathcal U}(\id\ot u^*)(y(u(f)))\,du,
 \quad\forall\; f\in F_n.$$
Similarly, we define $\wt z\,:\, F_n\to L_2^r(w_r;\el_2^n)$
associated to $z$. Let $\wt x=(\wt y,\,\wt z)$. By homogeneity,
$\|\wt y\|_{cb}\le c\|y\|_{cb}$ and $\|\wt z\|_{cb}\le
c\|z\|_{cb}$. On the other hand, both $\wt y$ and $\wt z$ commute
with all unitaries, so they must be multiples of $\id_{\el_2^n}$.
Therefore, there exist $a\in L_2(w_c)$ and $b\in L_2(w_r)$ such
that
 $$\wt y=a\ot \id\quad\mbox{and}\quad \wt z=b\ot\id\,.$$
Recall that $F_n$ is the subspace of constant functions of
$L_2^c(w_c;\el_2^n)+L_2^c(w_c;\el_2^n)$, which is just $\el_2^n$
at the algebraic level. Thus the inclusion map $\iota$ commutes
with all $u\in\mathcal U$. On the other hand, it is clear that the
quotient map $q$ also commutes with every $u$. Since $qx=\iota$,
we then deduce $q\wt x=\iota$. This amounts to saying that $a+b=1$
a.e. and concludes the proof.\cqd

\medskip

We are now ready to prove Theorem~\ref{ex}.

\medskip \n\emph{Proof of Theorem~\ref{ex}}. We keep the notations
in the proof of Theorem~\ref{c1s id explicit} and write
$F=K_{w_c,\, w_r}$ with the pair $(w_c,\, w_r)$ associated to the
fundamental functions of $F$. By Lemmas~\ref{exK} and \ref{diag},
we find
 $$ex(F_n)\sim \inf\left\{\|x\|_{cb}\,:\, x=(a\ot\id,\,b\ot\id),\;
 a\in L_2(w_c),\; b\in L_2(w_r),\; a+b=1\; {\rm a.e.}\right\}.$$
Fix an $x=(a\ot\id,\,b\ot\id)$ as above. Then
 $$\|x\|_{cb}\sim
 \|a\ot\id\|_{CB(F_n\,,\; L_2^c(w_c;\el_2^n))}\,+
  \|b\ot\id\|_{CB(F_n\,,\; L_2^r(w_r;\el_2^n))}\,.$$
It is easy to determine the two norms on the right hand side.
Indeed,   we have
 \be
 CB(F_n\,,\; L_2^c(w_c;\el_2^n))
 &=& L_2^c(w_c;\el_2^n)\ot_{\min}F_n^*\\
 &=& L_2^c(w_c)\ot_{\min}C^n\ot_{\min}F_n^*
 =L_2^c(w_c)\ot_{\min}C^n[F_n^*].
 \ee
It then follows that
 $$\|a\ot\id\|_{CB(F_n\,,\; L_2^c(w_c;\el_2^n))}
 =\F_{c, F^*}(n)^{1/2}\,\|a\|_{L_2(w_c)}\,.$$
Similarly,
  $$\|b\ot\id\|_{CB(F_n\,,\; L_2^r(w_r;\el_2^n))}
 =\F_{r, F^*}(n)^{1/2}\,\|b\|_{L_2(w_r)}\,.$$
Therefore, we deduce
 \be
 ex(F_n)
 &\sim&\inf\big\{
 \F_{c, F^*}(n)^{1/2}\,\|a\|_{L_2(w_c)}+
 \F_{r, F^*}(n)^{1/2}\,\|b\|_{L_2(w_r)}\,:\,
 a+b=1\; {\rm a.e.}\big\}\\
 &\sim& \big\|1\big\|_{\F_{c, F^*}(n)^{1/2}\,L_2(w_c)
 +\F_{r, F^*}(n)^{1/2}\,L_2(w_r)}\\
 &\sim& \big\|1\big\|_{L_2(\min(\F_{c, F^*}(n)w_c\,,\,
 \F_{r, F^*}(n)w_r))}\,.
 \ee
Now we apply the arguments from the proof of Theorem~\ref{c1s id
explicit}, where we also had to calculate the $L_2$-norm of $1$
with respect to the minimum of two weights. We outline the main
ingredients and leave the details to the reader. Indeed, on
$\rz_+$ we have $w_c=1$ and hence
 $$\min(\F_{c, F^*}(n)w_c\,,\,\F_{r, F^*}(n)w_r)=
 \min(\F_{c, F^*}(n)\,,\,\F_{r, F^*}(n)w_r).$$
The breaking point is given by $s_n\in\rz_+$ such that
 $$w_r(s_n)=\frac{\F_{c, F^*}(n)}{\F_{r, F^*}(n)}\,.$$
Then
 $$s_n\sim
 \F_{r, F}\big(\frac{\F_{r, F^*}(n)}{\F_{c, F^*}(n)}\big)$$
and by \eqref{g-h} and  \eqref{g-w}
 \be
 \big\|1\big\|^2_{L_2(\rz_+,\,\min(\F_{c, F^*}(n)w_c\,,\,
 \F_{r, F^*}(n)w_r))}
 &=&\int_0^{s_n}\F_{c, F^*}(n)ds + \int_{s_n}^\8 \F_{r, F^*}(n)w_r(s)ds\\
 &\sim& s_n\F_{c, F^*}(n)
 \sim \F_{r, F}\big(\frac{\F_{r, F^*}(n)}{\F_{c,
 F^*}(n)}\big)\F_{c, F^*}(n).
 \ee
Similarly, we have
 $$
 \big\|1\big\|^2_{L_2(\rz_-,\,\min(\F_{c, F^*}(n)w_c\,,\,
 \F_{r, F^*}(n)w_r))}
 \sim \F_{r, F^*}(n)\F_{c, F}\big(\frac{\F_{c, F^*}(n)}{\F_{r,
 F^*}(n)}\big).
 $$
Combining the preceding estimates, we obtain
 $$ex(F_n)^2\sim \F_{c, F^*}(n)\F_{r, F}\big(\frac{\F_{r, F^*}(n)}{\F_{c,
 F^*}(n)}\big) + \F_{r, F^*}(n)\F_{c, F}\big(\frac{\F_{c, F^*}(n)}{\F_{r,
 F^*}(n)}\big).$$
Together with Theorem~\ref{dual}, this proves the theorem. \cqd


\section{Examples}
 \label{Examples of applications}


We now apply the results in the previous sections to the column
and row $p$-spaces. It is known that $C_p\in QS(C\op R)$ for every
$1< p<\8$ (see \cite{ju-OH} and \cite{xu-embed}), so $C_p$ can be
represented as $K_{w_c,\,w_r}$ for a pair of weights on $\rz$ in
view of Corollary~\ref{hom 3 spaces}. In fact, it is easy to
calculate the fundamental functions of $C_p$. Thus
Theorem~\ref{uni} allows us to find a concrete representation of
$C_p$, i.e., to know explicitly the two weights $w_c$ and $w_r$.
Note that such concrete representations for $C_p$ are not new and
can be obtained by real interpolation (see \cite{xu-embed}). We
should also point out that a concrete representation of $OH$ was
first constructed  in \cite{ju-OH} by using Pusz-Woronowic formula
on the square root of a positive sesquilinear form on a Hilbert
space.

\medskip

We will be also interested in the sum and intersection of $C_p$
and $R_p$. Let
 $$CR_p=C_p+R_p \textrm{ for } 1\le p<2\quad\textrm{and}\quad
 CR_p=C_p\cap R_p \textrm{ for } 2\le p\le\8.$$
These spaces are still homogenous and Hilbertian, so belong to
$HQS(C\op R)$ too. Their importance stems from their links to
noncommutative Khintchine inequalities (see \cite{pis-ast}).

We start by calculating the fundamental functions of $C_p$ and
$CR_p$. The following  result is entirely elementary (see also
\cite[Lemma~5.9]{xu-embed}).

\begin{prop}\label{fund Cp}
 Let $1\le p\le\8$ and $p'$ be the conjugate index of $p$. Then
 $$\F_{c, C_p}(n)=n^{1/p'}\,,\quad \F_{r, C_p}(n)=n^{1/p}
 \quad\mbox{and}\quad
 \F_{c, CR_p}(n)=\F_{r, CR_p}(n)=n^{1/p'}\,.$$
 \end{prop}

\pf We have
 $$C[C_p]=(C[C],\; C[R])_{1/p}
 =(S_2,\; \mathbb K(\el_2))_{1/p}=S_{2p'}\,.$$
This implies $\F_{c, C_p}(n)=n^{1/p'}\,$. The second formula for
$C_p$ is obtained in the same way. On the other hand, for $1\le p<
2$ we have
  $$C[CR_p]=C[C_p]+ C[R_p]=S_{2p'}+S_{2p}=S_{2p'}\,;$$
whence $\F_{c, CR_p}(n)=n^{1/p'}\,$. The remaining formulas are
proved similarly. \cqd

\medskip

The result above shows that $C_p$ and $CR_p$ are regular for
$1<p<\8$. Thus all results in the previous sections apply to these
spaces. In particular, we can determine the completely $1$-summing
maps between them, their injectivity and exactness constants. We
collect all these in the following three statements. In the sequel
$\psi$ will denote the Orlicz function defined by $\psi(0)=0$ and
$\psi(t)=t^2\log(t+1/t)$ for $t>0$. It is clear that $\psi$
satisfies the $\Delta_2$-condition. The associated Orlicz space
$\el_\psi$ is traditionally denoted by $\el^2\log\el$. Note that
the inverse function of $\psi$ satisfies
 $$\psi^{-1}(t)\sim
\sqrt{2t}\,\big(\log\frac1t\big)^{-1/2} \quad\mbox{as}\quad
 t\to0.$$
It follows that the fundamental sequence $(\psi_n)$ of $\psi$ is
given by
 $$\psi_n \sim \sqrt{n\log(n+1)}\quad\mbox{as}\quad
 n\to\8.$$

\begin{thm}\label{c1c Cp}
 Let $1<p, q<\8$.
 \begin{enumerate}[\rm(i)]
 \item Let $r$ be determined by $2/r=1/p+1/q$.
 Then
  \be
  \Pi_1^o(C_p\,,\,C_{p'})=S_\psi
 \quad\textrm{and}\quad
  \Pi_1^o(C_p\,,\,C_q)=S_{\min(r,r')}\textrm{ for } q\neq p'\,.
  \ee
 \item Let $s$ be determined by $2/s=1/p+1/q'$. Then
  \be
  \Pi_1^o(CR_p\,,\,CR_{p})=S_\psi
 \quad\textrm{and}\quad
 \Pi_1^o(CR_p\,,\,CR_q)=S_{\min(2, s)}\textrm{ for } q\neq p\,.
\ee
 \end{enumerate}
 \indent Moreover, all relevant constants depend only on $p$ and
$q$.
 \end{thm}

\pf We consider only the part concerning $C_p$, the one on $CR_p$
being dealt with similarly.  By Theorem~\ref{c1s orlicz}, we know
that $\Pi_1^o(C_p\,,\,C_q)=S_\f$ for some Orlicz function $\f$.
Since $S_\f$ is completely determined by the fundamental sequence
of $\f$ (see Corollary~\ref{c1s funda}), we are reduced to
determine $\pi_1^o(\id_n\,:\, C_p^n\to C_q^n)$ for all $n$. This
is just simple integral calculations with help of Theorem~\ref{c1s
id explicit} and Proposition~\ref{fund Cp}. Indeed, we have
 \be
 &&\pi_1^o(\id_n\,:\, C_p^n\to C_q^n)^2
 \sim
 n^{1/p+1/q} + n^{1/p'+1/q'} +\\
 &&~~ n\Big[\int_1^{n^{1/p}} t^{-p/q'}dt +
 \int_1^{n^{1/p'}} t^{-p'/q}dt
 +\int_1^{n^{1/q'}} t^{-q'/p}dt +\int_1^{n^{1/q}} t^{-q/p'}dt\Big].
 \ee
If $q=p'$, then we deduce
 \be
 \pi_1^o(\id_n\,:\, C_p^n\to C_q^n)
 \sim \sqrt n + \sqrt {n\log(n+1)}\sim \sqrt {n\log(n+1)}\,.
 \ee
Since $( \sqrt {n\log(n+1)}\,)_{n\ge1}$ is the fundamental
sequence of $\psi$, we get $\Pi_1^o(C_p\,,\,C_{p'})=S_\psi$.

Assume $q>p'$. Then
 \be
 &&\pi_1^o(\id_n\,:\, C_p^n\to C_q^n)
 \sim
 n^{1/r} + n^{1/r'} +\\
 &&~~ n^{1/2}\Big[\int_1^{n^{1/p'}} t^{-p'/q}dt
 +\int_1^{n^{1/q'}} t^{-q'/p}dt\Big]^{1/2}
 \sim
 n^{1/\min(r, r')}\,.
 \ee
Again we are done. The case $q<p'$ is treated by symmetry on $p'$
and $q$. \cqd

\begin{rk}
 Tracking back the origin of the equivalence constants in
Theorem~\ref{c1c Cp}, one can find an explicit estimate for them
in terms of $p$ and $q$ and then one realizes that the result for
$q=p'$ may be obtained from that for $q\neq p'$ by a limit
procedure as $q\to p'$.
 \end{rk}

\begin{rk}
 It is easy to prove the inclusions
$\Pi_1^o(C_p,\, C_q)\subset S_{\min(r,r')}$ for $q\neq p'$ and
$\Pi_1^o(CR_p\,,\,CR_q)=S_{\min(2, s)}$ for $q\neq p$ are
contractive. Let us show the first one.
  \end{rk}

Let $u\in \Pi_1^o(C_p,\, C_q)$. Let $v: C\to C_p$ and $w: C_q\to
C$ be two finite rank maps. By the ideal property of completely
$1$-summing norms, we have
 $\pi_1^o(wuv)\le \|w\|_{cb}\,\pi_1^o(u)\,\|v\|_{cb}\,.$
However,
 $CB(C,\, C_p)=S_{2p}$ and $CB(C_q,\, C)=S_{2q}$.
It then follows that
 $$\|wuv\|_{1}\le \|w\|_{2q}\,\pi_1^o(u)\,\|v\|_{2p}\,.$$
Taking the supremum over all $v$ and $w$ in the unit balls of
$S_{2p}$ and $S_{2q}$, respectively, we obtain
 $$\|u\|_{r'}\le \,\pi_1^o(u)\,;$$
whence
 $$\Pi_1^o(C_p,\, C_q)\subset S_{r'}\,.$$
Note that the argument above remains valid with column spaces
replaced by row spaces. Thus
 $\Pi_1^o(R_{p},\,R_{q})\subset S_{r'}$ contractively.
However, $C_p=R_{p'}$ and $C_q=R_{q'}$. It then follows that
 $$\Pi_1^o(C_{p},\,C_{q})=\Pi_1^o(R_{p'},\,R_{q'})\subset S_{r}\,.$$
Therefore, we deduce
 $\Pi_1^o(C_p,\, C_q)\subset S_r\cap S_{r'}=S_{\min(r, r')}$
as desired.

\medskip

The following gives the injectivity constants of $C_p^n$ and
$CR_p^n$ for $p\neq2$. Recall that $C_2^n\cong CR_2^n\cong OH^n$
completely isometrically. Combining \cite{ju-OH} and
\cite{pisshlyak} we find
 $$\l_{cb}(C_2^n)=\l_{cb}(CR_2^n)
 =\l_{cb}(OH^n)\sim\frac{\sqrt n}{\sqrt{\log(n+1)}}\,.$$

\begin{thm}\label{inj Cp}
 Let $1< p<\8$ such that $p\neq2$.
Then
 $$\l_{cb}(C_p^n)\sim n^{\frac1{\max(p, \;p')}}
 \quad\mbox{and}\quad
 \l_{cb}(CR_p^n)\sim \frac{\sqrt n}{\sqrt{\log(n+1)}}$$
with equivalence constants depending only on $p$.
 \end{thm}

\pf This is an immediate consequence of Theorem~\ref{c1c Cp} and
Lemma~\ref{pi1-gama8}. Alternately, we can directly apply
Theorem~\ref{inj}. \cqd

\medskip

The second formula in the following is already in \cite{ju-hab}
(see Proposition~3.3.1.5 and Corollary~3.3.1.16 there), the first
has been known only for $p\in \{1,2,\infty\}$.

\begin{thm}\label{ex Cp}
 Let $1<p<\8$.
Then
 $$ex(C_p^n)\sim n^{\frac1{pp'}}
 \quad\mbox{and}\quad
 ex(CR_p^n)\sim  n^{\frac1{2p}}$$
with equivalence constants depending only on $p$.
 \end{thm}

\pf This follows immediately from Theorem~\ref{ex} and
Proposition~\ref{fund Cp}. \cqd

\begin{rk}
 It is clear that Theorem~\ref{ex Cp} holds for $p=1$ and $p=\8$.
It is also obvious that the first estimate in Theorem~\ref{inj Cp}
remains true for these endpoints. On the other hand, one has
$\l_{cb}(CR_1^n)=\l_{cb}(CR_\8^n)=\sqrt n$. Regarding
Theorem~\ref{c1c Cp}, one can check, without difficulty, that both
second equivalences in (i) and (ii) there hold if one of $p$ and
$q$ is $1$ or $\8$.
 \end{rk}

\bigskip

\n{\bf Acknowledgements.} This project got started in BIRS at
Banff in May/June 2004 while the authors were carrying out a
research team program there. They would like to thank BIRS for
providing excellent research facilities. They are also very
grateful to the anonymous referee for a careful reading of the
manuscript and for many helpful suggestions.

\bigskip


\bigskip

 \footnotesize{

\noindent M.J.: Department of Mathematics,  University of
Illinois, Urbana,
 IL 61801, USA\\
 junge@math.uiuc.edu\\

\noindent Q.X.: Laboratoire de Math{\'e}matiques, Universit{\'e} de
Franche-Comt{\'e},
25030 Besan\c{c}on Cedex,  France\\
qxu@univ-fcomte.fr\\}

\end{document}